\documentclass[11pt]{article}

\def\mathbb{\Bbb}
\usepackage{amsfonts,latexsym,amstext}
\usepackage{amsmath}
\usepackage{amssymb}

\usepackage{comment}

\newtheorem{theorem}{Theorem}[section]
\newtheorem{lemma}[theorem]{Lemma}
\newtheorem{proposition}[theorem]{Proposition}

\newtheorem{hypothesis}[theorem]{Hypothesis}

\newtheorem{remark}[theorem]{Remark}

\makeatletter
\@addtoreset{equation}{section}

\makeatother

\def\qed{{\hfill\hbox{\enspace${ \square}$}} \smallskip}
\def\sqr#1#2{{\vcenter{\vbox{\hrule height .#2pt \hbox{\vrule
 width .#2pt height#1pt \kern#1pt \vrule
width .#2pt} \hrule height .#2pt}}}}
\def\square{\mathchoice\sqr54\sqr54\sqr{4.1}3\sqr{3.5}3}

\def\ds{\begin{displaystyle}}
\def\eds{\end{displaystyle}}
\def\dis{\displaystyle }
\def\<{\langle }
\def\>{\rangle }

\def\R{\mathbb R}

\def\E{\mathbb E}
\def\P{\mathbb P}

\def\cala{{\cal A}}
\def\calb{{\cal B}}

\def\cale{{\cal E}}
\def\calf{{\cal F}}

\def\calk{{\cal K}}

\def\caln{{\cal N}}

\def\calp{{\cal P}}

\def\calu{{\cal U}}

\def\call{{\cal L}}

\def\K{\mathbb K}
\def\H{\mathbb H}

\title{Backward stochastic differential equations
and optimal control of marked point processes}

\date{}
\author{Fulvia Confortola, Marco Fuhrman
\\Politecnico di Milano,
Dipartimento di Matematica\\
piazza Leonardo da Vinci 32, 20133 Milano, Italy\\
e-mail: fulvia.confortola@polimi.it, marco.fuhrman@polimi.it}

\begin{document}

\maketitle

\begin{abstract}
We study a class of backward
stochastic differential equations (BSDEs) driven by a random measure or,
equivalently, by a marked point process. Under appropriate assumptions
we prove well-posedness and continuous dependence of the solution
on the data.
We next address  optimal control problems for point processes of general non-markovian type
and show  that  BSDEs can be used to prove existence of an optimal control
and to represent the value function. Finally we introduce a Hamilton-Jacobi-Bellman equation, also
stochastic and of backward type,
for this class of control problems: when the state  space is finite or countable we show that it
 admits a unique solution which identifies the (random) value function and can be
represented by means of the BSDEs introduced above.
\end{abstract}



\section{Introduction}

The purpose of this paper is to study a class of backward
stochastic differential equations (BSDEs for short) and apply
these results to solve optimal control problems for
marked point processes. Under appropriate assumptions, an associated
Hamilton-Jacobi-Bellman equation of stochastic type is also
introduced and solved in this non-markovian framework.

General nonlinear BSDEs driven by the Wiener process were first
solved in \cite{PaPe}. Since then, many generalizations have been
considered where the Wiener process was replaced by more
general processes. Among the earliest results we mention in
particular \cite{EPQ}, \cite{ElKHua}, to which some of our results
are inspired,  and we refer e.g. to \cite{CaFeSa}
 for a recent result and for indications on the
existing bibliography.

We address a class of BSDEs  driven by a random measure, naturally
associated to a marked point process. There exists a large
literature on this class of processes, and in particular to the corresponding
optimal control problems: we only mention the
classical treatise \cite{B} and the recent book \cite{BrLa} as
general references. In spite of that, there are relatively few
results on their connections with BSDEs. In the general
formulation of a BSDE driven by a random measure, one of the
unknonwn processes (the one associated with the martingale part,
or $Z$-process) is in fact a random field. This kind of equations
has been introduced in \cite{TaLi},
and has been later considered in \cite{BaBuPa}, \cite{Roy} in the markovian case,
where the associated (nonlocal) partial differential equation and related
non-linear expectations
have
been studied.

 In these papers the BSDE contains a diffusive part and
a jump part,
but the latter is only considered in the case of a Poisson
random measure.  In order to give a probabilistic representation
of solutions to quasi-variational inequalities in the theory of stochastic impulse control,
in \cite{KhMaPhZh}  a more difficult problem involving also constraints
on the jump part is formulated and solved, but still in the Poisson case and
in a markovian framework.

To our knowledge, the only general
result beyond the Poisson case  is the paper \cite{Xia}.
Here, under conditions of Lipschitz type on the coefficients and assuming
the validity of appropriate martingale representation theorems,
a general  BSDE driven by  a diffusive and a jump part is considered and
well-posedness results and a comparison theorem are proven.
However, it seems that in this paper the formulation of the BSDE was not chosen in view of
applications to optimal control problems. Indeed, in contrast to
\cite{TaLi} or
\cite{BaBuPa},
the generator of the BSDE depends on the $Z$-process in a specific
way (namely as an integral of a Nemytskii operator)
that is generally not valid for the hamiltonian function of  optimal control problems
(compare for instance formula  \eqref{defhamiltonianuno} below) and therefore
prevents direct applications to these problems.

In our paper  we consider a BSDE driven by a random measure,
without diffusion part, on a finite time interval, of the following form:
\begin{equation}\label{BSDEuno}
Y_t + \int_t^T \int_KZ_s(y)\,q(ds \,dy)= \xi + \int_t^T
f_s(Y_s,Z_s(\cdot))\, dA_s, \qquad t\in  [0,T],
\end{equation}
where the generator $f$ and the final condition $\xi$ are given.

Here the basic probabilistic datum is a marked point process $(T_n,\xi_n)$ where
$(T_n)$ is an increasing sequence of random times and $(\xi_n)$ a sequence of random
variables in the state (or  mark) space $K$. The corresponding random counting measure
is $p(dt\,dy)= \sum_{n} \delta_{(T_n,\xi_n)} $, where $\delta$ denotes the Dirac measure.
We denote $(A_t)$ the compensator of the counting process $(p([0,t]\times K))$ and by
$\phi_t(dy)\,dA_t$ the (random) compensator of  $p$. Finally, the compensated measure
$q(dt\, dy)= p(dt\,dy) -
 \phi_t(dy)\,dA_t$ occurs in equation \eqref{BSDEuno}.
The unknown process is a pair
 $(Y_t,Z_t(\cdot))$, where $Y$ is a real progressive process and
$\{Z_t(y),\, t\in [0,T],\,y\in K\}$ is a predictable random field.

The random measure $p$ is fairly general, the
only restriction  being  non explosion (i.e. $T_n\to\infty$) and
the requirement that $(A_t)$ has continuous trajectories.
We allow
the  space $K$ to be of general type, for instance a Lusin space.
Therefore our results can also be directly applied to marked point processes with discrete state
space. We mention at this point that the specific case  of finite or countable Markov chains
has been studied in \cite{Coh-Ell-1}, \cite{Coh-Ell-2}, see also
\cite{Coh-Ell-3}
for generalizations.

The basic hypothesis on the generator $f$ is a Lipschitz
condition  requiring that for some constants $L \ge 0$, $L' \ge 0$,
$$
|f_t(\omega,r,z(\cdot))-f_t(\omega,r',z'(\cdot))|\leq L'|r-r'|
+ L\left(\int_K |z(y)-z'(y)|^2\phi_t(\omega,dy)\right)^{1/2}
$$
for all $(\omega,t)$, for $r,r'\in\R$, and $z,z'$
in appropriate function spaces  (depending on $(\omega,t)$):  see below for precise
statements.
  We note that the generator of the BSDE can depend on the unknown
$Z$-process in a general functional way: this is required in the applications to
optimal control problems that follow, and it is shown that our assumptions can be
 effectively verified in a number of cases. In order to solve the equation, beside measurability assumptions,
  we require
 the summability condition
 $$
 \E \int_0^T e^{\beta A_t}|f_t(0,0)|^2 dA_t+
 \E \, [e^{\beta A_T}|\xi|^2] < \infty,
 $$
 to hold for some $\beta>L^2+2L'$. Note that in the Poisson case mentioned above we have
 a deterministic compensator $\phi_t(dy)\,dA_t=\pi(dy)\,dt$ for some fixed measure $\pi$
 on $K$  and the summability condition reduces to a simpler form, not involving exponentials
 of stochastic processes.
We
prove existence, uniqueness, a priori estimates and continuous dependence upon the
data for the solution to the BSDE.

The results described so far are presented
in section
\ref{sec-backward}, after an introductory  section
 devoted
to notation and preliminaries.

In section \ref{sec-control} we formulate a class of optimal
control problems for marked point processes, following a classical
approach exposed for instance in \cite{B}.
For every fixed $(t,x)\in [0,T]\times K$,
the cost to be minimized and the corresponding   value function
are
$$
J_t(x,u(\cdot))=\E_u^{\calf_t}\left[
\int_t^Tl_s(X_s^{t,x},u_s)\;dA_s + g(X_T^{t,x})\right],
\qquad v(t,x)= \mathop{\rm ess\ inf}_{u(\cdot)\in\cala }J_t(x,u(\cdot)),
$$
where $\E_u^{\calf_t}$ denotes the conditional expectation with respect to a new probability $\P_u$,
depending on a control process $(u_t)$ and
defined by means of an absolutely continuous change of
measure: the choice of the
 control process  modifies the
compensator of the random measure under $\P_u$ making it equal to
$ r_t (y,u_t)\phi_t(dy)\,dA_t$ for some given function $r$. To this control problem
we associate the BSDE
\begin{equation}\label{bsdecontrollouno}
    Y_s^{t,x}+\int_s^T\int_KZ_r^{t,x}(y)\,q(dr\,dy) =
g(X_T^{t,x}) +\int_s^T f(r,X_r^{t,x},Z_r^{t,x}(\cdot))\,dA_r,
\qquad s\in [t,T].
\end{equation}
where $(X_r^{t,x})$ is a family of marked point processes, each starting from $x$
at time $t$, and
the  generator contains the hamiltonian function
\begin{equation}\label{defhamiltonianuno}
    f(\omega,t,x,z(\cdot))=\inf_{u\in U}\left\{
l_t(\omega, x,u)+ \int_K z(y) \, (r_t
(\omega,y,u)-1)\,\phi_t(\omega,dy)\right\}.
\end{equation}
Assuming that the infimum is in fact a minimum, admitting a suitable selector, together
with a summability condition of the form
$$
\E \exp \left(\beta A_T\right)+
 \E[|g(X_T^{t,x})|^2e^{\beta A_T} ]<\infty
$$
for a sufficiently large value of $ \beta$,
we prove that the optimal control
problem has a solution, and that the  value function and the optimal
control can be represented by means of the solution to the BSDE.

We note that optimal control of point processes is a classical topic in stochastic analysis,
and the first main contributions date back several decades: we refer the reader for instance to
the corresponding chapters of the treatises \cite{B} and  \cite{E}. The markovian case
has been further investigated in depth, even for more general classes of processes,
see e.g. \cite{Da-bo}. The results we present in this paper are an attempt toward
an alternative systematic approach, based on BSDEs. We hope this may lead to useful results in
the future, for instance in connection with computational issues and a
better understanding of the nonmarkovian situation.
Although this approach is analogous to the diffusive case, it
seems that it is pursued here for the first time in the case of marked point
processes.   In particular it
differs from the control-theoretic applications
addressed in  \cite{TaLi}, devoted to a version of the stochastic maximum principle.

Finally, in section \ref{sec-HJB}, we introduce the following
Hamilton-Jacobi-Bellman equation (HJB for short) associated to the optimal
control problem described above:
\begin{equation}\label{HJBstochuno}
\begin{array}{l}
  \dis
  v(t,x)+ \int_t^T\int_KV(s,x,y)\,q(ds\, dy)
  \\\qquad\dis
  =
  g(x)
  +
  \int_t^T
\int_K\Big( v(s,y)-v(s,x) + V(s,y,y)-V(s,x,y)\Big)\,\phi_s(dy) \,dA_s
\\\qquad\dis
+\int_t^T f\Big(s,x,v(s,\cdot)-v(s,x)+ V(s,\cdot,\cdot)\Big)\,dA_s,
\qquad t\in [0,T],\, x\in K,
\end{array}
\end{equation}
where $f$ be the hamiltonian function   defined in
(\ref{defhamiltonianuno}). The solution is a pair of random fields
$\{v(t,x),V(t,x,y)\,:\, t\in [0,T],\, x,y\in K\}$, and
in this non-markovian framework the HJB equation
is stochastic and of backward type, driven by the same random measure as before. Thus,
the previous results are applied to prove its well-posedness. For
technical reasons, however, we limit ourselves to the case where
the state space $K$  is at most countable: although this is
a considerable restriction with respect to the previous results,
it allows to treat important classes of control problems, for instance
those related to queuing systems. Under appropriate assumptions, similar
to those outlined above, we prove
that the HJB equation is well-posed and that $v(t,x)$ coincide with the
(stochastic) value function of the optimal control problem and it can be
represented by means of the associated BSDE. A backward stochastic
HJB equation has been first introduced in \cite{Pe} in the
diffusive case, where the corresponding theory is still not
complete due to greater technical difficulties.  It is an interesting fact that
the parallel case of jump processes can be treated using BSDEs and fairly complete
results can be given, at least under the restriction mentioned above:
this is perhaps due to the different nature of the control problem (here
the laws of the controlled processes are obtained via an absolutely continuous
change of measure,  in  contrast to \cite{Pe}).
We borrow some ideas from \cite{Pe}, in particular the use of a
formula of Ito-Kunita type proven below, that suggested the unusual
form of \eqref{HJBstochuno}.
We are not aware
of any previous result on backward HJB equations in a
non-diffusive context.

The results of this paper admit several variants and
generalizations: some of them are not included here for reasons of brevity and some
are presently in
preparation. For instance, the BSDE approach to optimal control of Markov jump
processes deserves a specific treatment; moreover,
  BSDEs driven by random measures can be studied without
Lipschitz assumptions on the generator, along the lines of the many results available
in the diffusive case, or extensions to the case of vector-valued
process $Y$ or of random time interval can be considered.

\section{Notations, preliminaries and basic assumptions}
\label{notations}

In this section we are going to recall basic notions
on marked points processes, random measures and corresponding stochastic
integrals, that will be constantly used in the rest of the paper.
We also formulate several assumptions that will remain in force throughout.

\subsection{Marked point processes}
\label{sectionmpp}

Let $(\Omega,\calf, \P)$ be a complete probability space
and $(K,\calk)$ a measurable space.
Assume we have a sequence $(T_n,\xi_n)_{n\ge 1}$ of
random variables, $T_n$ taking values in $[0,\infty]$
and $\xi_n$ in $K$. We set $T_0=0$
and we assume, $\P$-a.s.,
$$T_n<\infty\; \Longrightarrow\; T_n<T_{n+1},\qquad n\ge 0.
$$
We call $(T_n)$ a point process and $(T_n,\xi_n)$ a marked
point process. $K$ is called the mark space, or state space.

In this paper we will always assume that $(T_n)$  is nonexplosive,
i.e.
$T_n\to\infty$ $\P$-a.s.

For every $A\in\calk$ we define the counting processes
$$
N_t(A)= \sum_{n\ge 1}1_{T_n\le t}1_{\xi_n\in A}, \qquad t\ge 0
$$
and we set $N_t=N_t(K)$. We define the filtration generated
by the counting processes by first introducing the $\sigma$-algebras
$$
\calf^0_t=\sigma(N_s(A)\;:\;s\in [0,t],\, A\in \calk), \qquad t\ge 0,
$$
and setting
$$
\calf_t=\sigma(\calf_t^0, \caln), \qquad t\ge 0,
$$
where  $\caln$ denotes the family of $\P$-null sets in $\calf$.
It turns out that $(\calf_t)_{t\ge 0}$ is right-continuous and therefore
satisfies the usual conditions.
In the following all measurability concepts for stochastic processes (e.g. adaptedness,
predictability) will refer to
 the filtration $(\calf_t)_{t\ge 0}$.  The predictable $\sigma$-algebra
 (respectively, the progressive $\sigma$-algebra)
 on $\Omega\times [0,\infty)$  will be denoted by $\calp$ (respectively,
 by $Prog$). The same symbols will also denote the
 restriction to
 $\Omega\times [0,T]$ for some $T>0$.

It is known that there exists an increasing, right-continuous predictable  process $A$
satisfying $A_0=0$ and
$$
\E\int_{0}^\infty H_t\; dN_t=\E\int_{0}^\infty H_t\; dA_t
$$
for every nonnegative predictable process $H$. The above
stochastic integrals are defined for $\P$-almost every $\omega$ as
ordinary (Stieltjes) integrals. $A$ is called the compensator, or
the dual predictable projection, of $N$. In the following we will
always make the basic assumption that $\P$-a.s.
\begin{equation}\label{acontinuo}
    A \;\text{has continuous trajectories}
\end{equation}
which are in particular finite-valued.

We finally fix $\xi_0\in K$ (deterministic)
and we define
\begin{equation}\label{defxdibase}
X_t=\sum_{n\ge 0}\xi_n\,1_{[T_n,T_{n+1})}(t), \qquad t\ge 0.
\end{equation}
We do not assume that $\P(\xi_n\neq \xi_{n+1})=1$. Therefore
in general trajectories of $(T_n,\xi_n)_{n\ge 0}$ cannot be reconstructed
from trajectories of $(X_t)_{t\ge 0}$ and the filtration $(\calf_t)_{t\ge 0}$
is not the natural completed filtration of  $(X_t)_{t\ge 0}$.

\subsection{Random measures and their compensators}

For $\omega\in\Omega$ we define a measure on $((0,\infty)\times K,
\calb((0,\infty))\otimes \calk)$ setting
$$
p(\omega,C)=\sum_{n\ge 1} 1_{(T_n(\omega),\xi_n(\omega))\in C}
\qquad C\in \calb((0,\infty))\otimes \calk,
$$
where $\calb(\Lambda)$ denotes   the Borel $\sigma$-algebra of any
topological space $\Lambda$.   $p$ is called a random measure since
 $\omega\mapsto p(\omega,C)$ is $\calf$-measurable for
fixed $C$. We also use the notation $p(\omega, dt\,dy)$ or
$p(dt\,dy)$. Notice that $p((0,t]\times A)=N_t(A)$ for $t> 0,
A\in\calk$.

Under mild assumptions on $K$  it can be proved that there exists a function
$\phi_t(\omega,A)$ such that

\begin{enumerate}
\item for every $\omega\in\Omega$, $t\in [0,\infty)$, the
mapping $A\mapsto \phi_t(\omega,A)$ is a probability measure on $(K,\calk)$;

\item for every $A\in\calk$, the process $(\omega,t)\mapsto \phi_t(\omega,A)$
is predictable;

\item for every nonnegative $H_t(\omega,y)$, $\calp\otimes \calk$-measurable,
we have
$$
\E \int_{0}^\infty H_t(y)\; p(dt\,dy)=\E \int_{0}^\infty H_t(y)\;\phi_t(dy)\, dA_t.
$$
\end{enumerate}

 For instance, this holds  if
 $(K,\calk)$ is a Lusin space
with its Borel $\sigma$-algebra (see \cite{J} Section 2),
but since the Lusin property will not play any further role below,
in the following we will simply assume the existence of
$\phi_t(dy)$ satisfying $1$-$2$-$3$ above.

The random measure $\phi_t(\omega,dy)\, dA_t(\omega)$ will be
denoted $\tilde{p}(\omega, dt\,dy)$, or simply
$\tilde{p}(dt\,dy)$, and will be called the compensator, or the
dual predictable projection, of $p$.

\subsection{Stochastic integrals}\label{subsec-stocint}

Fix  $T>0$,  and let $H_t(\omega,y)$ be a $\calp\otimes
\calk$-measurable real function satisfying
$$
\int_0^T \int_K |H_t(y)| \;\phi_t(dy)\, dA_t <\infty, \qquad
\P-a.s.
$$
Then the following stochastic integral can be defined
\begin{equation}\label{defstocint}
  \int_0^t \int_K H_s(y)\; q(ds\,dy):= \int_0^t \int_K H_s(y)\;
p(ds\,dy)-\int_0^t \int_K H_s(y) \;\phi_s(dy)\, dA_s , \qquad t
\in [0,T],
\end{equation}
   as the difference of ordinary
integrals with respect to $ p$ and $\tilde{p}$. Here and in the
following the symbol $\int_a^b$ is to be understood as an integral
over the interval ${(a,b]}$. We shorten this identity writing
$q(dt \, dy)= p(dt\, dy) -\tilde{p}(dt\, dy)= p(dt\,
dy)-\phi_t(dy)\, dA_t$.
 Note
that
$$
\int_0^t \int_K H_s(y)\; p(ds\,dy) = \sum_{n\ge 1, T_n\le t}
H_{T_n}(\xi_n)
$$
is always well defined since we are assuming that $T_n\to\infty$
$\P$-a.s.

 For $r\ge 1$ we define $\call^{r,0}(p)$ as
the space of $\calp\otimes \calk$-measurable real functions
$H_t(\omega,y)$ such that
$$
\E \int_{0}^T \int_K |H_t(y)|^r\; p(dt\,dy)=
\E \int_{0}^T \int_K  |H_t(y)|^r\;\phi_t(dy)\, dA_t<\infty
$$
(the equality of the integrals follows from the definition of
$\phi_t(dy)$). Given an element $H$ of $\call^{1,0}(p)$, the
stochastic integral \eqref{defstocint}
 turns out to be a  a finite variation martingale.

The key result used in the construction of a solution to the BSDE
\eqref{BSDE} is the integral representation theorem of marked
point process martingales (see e.g.  \cite{Da-art},\cite{Da-bo}),
which is a counterpart of the well known representation result for
Brownian martingales (see e.g. \cite{Re-Yor} Ch V.3  or  \cite{E}
Thm 12.33). Recall that $(\calf_t)$ is the filtration generated by
the jump process, augmented in the usual way.

\begin{theorem}\label{rappresentazione}
Let $M$ be a cadlag $(\calf_t)$-martingale on $[0,T]$. Then    we
have
$$M_t= M_0 + \int_0^t \int_K H_s(y)\, q(ds\; dy),\qquad   t \in [0, T],
$$
for some process $H \in \mathcal{L}^{1,0}(p)$.
\end{theorem}

\subsection{A family of marked point processes.}
\label{famigliaparametrizzata}

In the following, in order to use dynamic programming arguments, it will be useful
to introduce a family of processes instead of the single process $X$, each
starting at a different time from different points.

Let $(T_n,\xi_n)$ be the marked
point process introduced in section \ref{sectionmpp}.
We fix $t\ge 0$ and we introduce  counting processes relative to the
time interval $[t,\infty)$ setting
$$
N_s^t(A)= \sum_{n\ge 1}1_{t<T_n\le s}1_{\xi_n\in A}, \qquad s\in [t,\infty), \;A\in\calk,
$$
and   $N_s^t=N_s^t(K)$. Then $N_s^t(A)=p^t((t,s]\times A)$ for
$s>t, A\in\calk$, where the random measure $p^t$ is the
restriction of $p$ to  $(t,\infty)\times K$. With these
definitions it is easily verified that the compensator of $p^t$
(respectively, $N^t$) is the restriction of $\phi_s(dy)\,dA_s$
(respectively, $A$) to $ [t,\infty)\times K$ (respectively, $
[t,\infty)$).

Now we fix $t\ge 0$ and $x\in K$. Noting that $N_t$ is the number of jump times $T_n$
in the interval $[0,t]$, so that $T_{N_t}\le t<T_{N_t+1}$, we define
$$
X_s^{t,x}= x\,1_{[t,T_{N_t+1})}(s)+ \sum_{n\ge N_t+1} \xi_n\, 1_{[T_n,T_{n+1})}(s),
\qquad s\in [t,\infty).
$$
In particular, recalling the definition of
the process $X$, previously defined by formula (\ref{defxdibase}) and starting at
 point $\xi_0\in K$, we observe that
$X=X^{0,\xi_0}$.

For arbitrary $t,x$ we also have $X_s^{t,x}=X_s$ for $s\ge T_{N_t+1}$ and,
finally, for $0\le u\le t\le s$ and $x\in K$ the identity
$
X_s^{t,X_t^{u,x}}=X_s^{u,x}
$
is easy to verify.

\section{The backward equation} \label{sec-backward}

From now on, we fix a deterministic terminal time $T > 0$.

For given $\omega\in\Omega$ and $t\in[0,T]$, we denote $
\call^r(K,\mathcal{K}, \phi_t(\omega,dy) )$   the usual space of
$\calk$-measurable maps $z:K\to\R$ such that $\int_K
|z(y)|^r\phi_t(\omega,dy)<\infty$ (below we will only use $r=0$ or
$1$).

Next we introduce several classes of stochastic processes,
depending on a parameter $\beta >  0$.
\begin{itemize}
  \item
$\call_{Prog}^{2, \beta}(\Omega\times [0,T])$ denotes the set of
real progressive processes $Y$  such that
$$|Y|_{\beta}^2:=\E\int_0^T e^{\beta A_t }|Y_t|^2dA_t < \infty.
$$

\item $\mathcal{L}^{2,\beta}(p)$ denotes the set of mappings
$Z : \Omega\times [0, T] \times K \rightarrow \R$ which are
$\mathcal{P}\otimes \mathcal{K}$-measurable
and such that
$$ \|Z\|_{\beta}^2:= \E \int_0^T
\int_Ke^{\beta A_t } |Z_t(y)|^2 \phi_t(dy)dA_t < \infty.$$
\end{itemize}

We say that $Y,Y'\in \call_{Prog}^{2, \beta}(\Omega\times [0,T])$
(respectively, $Z,Z'\in\mathcal{L}^{2,\beta}(p)$) are equivalent
if
 they coincide almost everywhere with respect to the measure
$dA_t(\omega)\P(d\omega)$ (respectively, the measure
$\phi_t(\omega,dy)dA_t(\omega)\P(d\omega)$) and this happens if
and only if $|Y-Y'|_{\beta}=0$ (respectively,
$\|Z-Z'\|_{\beta}=0$). We denote $L_{Prog}^{2, \beta}(\Omega\times
[0,T])$ (respectively, ${L}^{2,\beta}(p)$) the corresponding set
of equivalence classes, endowed with the norm $|\cdot|_{\beta}$
(respectively, $\|\cdot\|_{\beta}$). $L_{Prog}^{2,
\beta}(\Omega\times [0,T])$ and ${L}^{2,\beta}(p)$ are Hilbert
spaces, isomorphic to $ L^{2, \beta}(\Omega\times [0,T], Prog,
e^{\beta A_t(\omega) }\,dA_t(\omega)\,\P(d\omega) )$ and $L^{2,
\beta}(\Omega\times [0,T]\times K, \calp\otimes\calk, e^{\beta
A_t(\omega) }\,\phi_t(\omega,dy)\,dA_t(\omega)\,\P(d\omega) ) $
respectively.

Finally we introduce the Hilbert space
 $\K^{\beta}= L_{Prog}^{2, \beta}(\Omega\times
[0,T]) \times {L}^{2, \beta}(p)$, endowed with the norm
$||(Y,Z)||^{2}_{ \beta}:=|Y|_{\beta}^2+||Z||_{\beta}^2$.

In the following we will consider the backward stochastic
differential equation: $\P$-a.s.,
\begin{equation}\label{BSDE}
Y_t + \int_t^T \int_KZ_s(y)\,q(ds \,dy)= \xi + \int_t^T
f_s(Y_s,Z_s(\cdot))\, dA_s, \qquad t\in  [0,T],
\end{equation}
where the generator $f$ and the final condition $\xi$ are given and and we look
for unknown processes $(Y,Z)\in\K^{\beta}$.

Let us consider the following assumptions on the data $f$ and
$\xi$.

\begin{hypothesis}\label{hyp:BSDE}
\begin{enumerate}
\item The final condition $\xi:\Omega\to\R$ is $\calf_T$-measurable and
$\E \, e^{\beta A_T}|\xi|^2 < \infty$.

\item For every $\omega \in \Omega$,
$t \in [0,T]$, $r \in \R$, a mapping  $f_t(\omega,r,
\cdot):\call^{2}(K, \calk,\phi_t(\omega,dy))\to \R$ is given,
satisfying the following assumptions:
\begin{itemize}
\item[(i)] for every $ Z \in \mathcal{L}^{2,\beta}(p)$
the mapping
\begin{equation}\label{misurabilitaf}
  (\omega,t,r)\mapsto f_t(\omega,r, Z_t(\omega,\cdot))
\end{equation}
is $Prog\otimes \calb(\R)$-measurable;
\item[(ii)] there exists $L \ge 0$, $L' \ge 0$
such that for every $\omega\in\Omega$, $t \in [0,T]$, $ r, r' \in
\R$, $z,z' \in \call^2(K,\mathcal{K}, \phi_t(\omega,dy) )$ we have
\begin{equation}\label{generatorelip}
|f_t(\omega,r,z(\cdot))-f_t(\omega,r',z'(\cdot))|\leq L'|r-r'|
+ L\left(\int_K |z(y)-z'(y)|^2\phi_t(\omega,dy)\right)^{1/2};
\end{equation}
\item[(iii)] We have
\begin{equation}\label{generatoresommabile}
\E \int_0^T e^{\beta A_t}|f_t(0,0)|^2 dA_t <\infty.
\end{equation}
\end{itemize}
\end{enumerate}
\end{hypothesis}

\begin{remark} \begin{em}
\begin{enumerate}
\item
The slightly involved measurability condition on the generator
seems unavoidable, since the mapping  $f_t(\omega,r, \cdot)$ has a
domain which depends on $(\omega,t )$. However, in the following
section, we will see how it can be effectively verified in
connection with   optimal control problems.

 Note that if $ Z \in
\mathcal{L}^{2,\beta}(p)$ then $ Z_t(\omega,\cdot)$ belongs to
$\call^2(K,\mathcal{K}, \phi_t(\omega,dy) )$ except possibly on a
predictable set of points $(\omega,t )$ of measure zero with
respect to $dA_t(\omega)\P(d\omega)$, so that the requirement on
the measurability of the map \eqref{misurabilitaf} is meaningful.

\item We note the inclusion
\begin{equation}\label{inclusioneprocessi}
   \mathcal{L}^{2,\beta}(p)\subset  \mathcal{L}^{1,0}(p).
\end{equation}
Indeed if $Z\in \mathcal{L}^{2,\beta}(p)$ then the inequality
$$
\int_0^T \int_K |Z_t(y)|\,\;\phi_t(dy)\, dA_t \le \left( \int_0^T
\int_K |Z_t(y)|^2\,\;\phi_t(dy)\,e^{\beta A_t} dA_t\right)^{1/2}
\left( \int_0^T e^{-\beta A_t} dA_t\right)^{1/2}
$$
and the fact that $\int_0^T e^{-\beta A_t}
dA_t=\beta^{-1}(1-e^{-\beta A_T})\le \beta^{-1}$ imply that $Z\in
\mathcal{L}^{1,0}(p)$.

 It
follows from \eqref{inclusioneprocessi} that the martingale $M_t=
\int_0^t \int_K Z_s(y)\, q(ds\; dy)$ is well defined for $Z\in
\mathcal{L}^{2,\beta}(p)$ and has cadlag trajectories $\P$-a.s. It
is easily checked that $M$ only depends on the equivalence class
of $Z$ as defined above.
\end{enumerate}
\end{em}
\end{remark}

\begin{lemma} \label{BSDEcasobanale}
Suppose that $f:\Omega\times [0,T]\to \R$ is  progressive,
$\xi:\Omega\to \R$ is $\calf_T$-measurable and
$$
\E\, e^{\beta A_T}|{\xi}|^2 + \E \int_0^T e^{\beta A_s}|{f}_s|^2
dA_s<\infty
$$
for some $\beta> 0$. Then there exists a unique pair $(Y,Z)$ in $\K^{\beta}$ solution to
the BSDE
\begin{equation}\label{BSDElineare}
  Y_t+\int_t^T \int_K Z_s(y) \, q(ds \;dy)=
   \xi +\int_t^Tf_s  \, dA_s.
\end{equation}
Moreover the following identity holds:
\begin{equation}\label{identitaenergia}
\begin{array}{l}\dis
\E\, e^{\beta A_t}| {Y}_t|^2 +\beta\,\E\int_t^T   e^{\beta
A_s}|{Y}_s|^2 dA_s
+ \E\int_t^T \int_K e^{\beta A_s}|{Z}_s(y)|^2 \phi_s(dy) dA_s \\
\dis\qquad =\E\, e^{\beta A_T}|{\xi}|^2
  +2\E \int_t^T  e^{\beta A_s} {Y}_{s} \,f_s dA_s,
  \qquad t\in [0,T],
\end{array}
\end{equation}
and and there exist
two constants $c_1(\beta)=4(1+\frac{1}{\beta})$ and $c_2(\beta)= \frac{8}{\beta}(1+\frac{1}{\beta})$ such that
\begin{equation}\label{stimaprlineare}
 \E\int_0^T   e^{\beta A_s}| {Y}_s|^2 dA_s
+ \E\int_0^T \int_K e^{\beta A_s}| {Z}_s(y)|^2 \phi_s(dy) dA_s \le c_1(\beta)
\E\,e^{\beta A_T}|\xi|^2 +  c_2(\beta)\int_0^Te^{\beta A_s}|f_s|^2  \, dA_s.
\end{equation}
\end{lemma}

\noindent{\bf Proof.}
Uniqueness follows immediately using the linearity of \eqref{BSDElineare} and taking the
conditional expectation given $\calf_t$.

Assuming that $(Y,Z)\in\K^{\beta}$ is a solution, we now prove the identity
\eqref{identitaenergia}.
From the Ito formula applied to $ e^{\beta A_t}| {Y}_t|^2$ it follows that
$$d( e^{\beta A_t}| {Y}_t|^2)= \beta  e^{\beta A_t}| {Y}_t|^2 dA_t + 2 e^{\beta A_t} {Y}_{t-} dY_t + e^{\beta A_t}|\Delta {Y}_t|^2.$$

So integrating on $[t,T]$ and recalling that $A$ is continuous,
\begin{eqnarray}
\nonumber e^{\beta A_t}| {Y}_t|^2 & = & - \int_t^T \beta  e^{\beta A_s}| {Y}_s|^2 dA_s - 2 \int_t^T e^{\beta A_s} {Y}_{s-} \int_K {Z}_s(y) q(ds\,dy) - \sum_{t<s\le T} e^{\beta A_s}|\Delta  {Y}_s|^2\\
\label{eq:int-Itobis} &  & +e^{\beta A_T}| {\xi}|^2
+2\int_t^T  e^{\beta A_s} {Y}_{s}  \,f_s\, dA_s.
\end{eqnarray}
The integral process $\int_0^t  e^{\beta A_s} {Y}_{s-} \int_K Z_s(y) q(ds\,dy)$
is a martingale, because the integrand process $ e^{\beta A_s} {Y}_{s-}   {Z}_s(y)$ is in $\call^1(p)$:
 in fact from the Young inequality we get
\begin{eqnarray*}
             \E \int_0^T \int_K  e^{\beta A_s}| {Y}_{s-}| |  {Z}_s (y)| \phi_s(dy)dA_s && \\
     \le \frac{1}{2}\E \int_0^T e^{\beta A_s}| {Y}_{s-}|^2 dA_s &+&
      \frac{1}{2}\E \int_0^T \int_K  e^{\beta A_s}|  {Z}_s (y)|^2 \phi_s(dy)dA_s < +\infty.
             \end{eqnarray*}
Moreover  we have
\begin{eqnarray*}
\sum_{0<s\le t} e^{\beta A_s}|\Delta  {Y}_s|^2 & = & \int_0^t \int_K e^{\beta A_s}| {Z}_s(y)|^2p(ds\,dy)
  \\
 & = &\int_0^t \int_K e^{\beta A_s}| {Z}_s(y)|^2 q(ds\,dy)+
  \int_0^t \int_K e^{\beta A_s}| {Z}_s(y)|^2 \phi_s(dy) dA_s,
 \end{eqnarray*}
where the stochastic integral with respect to $q$ is a martingale.
Taking the expectation in \eqref{eq:int-Ito} we obtain
\eqref{identitaenergia}.

We now pass to the proof of existence of the required solution.
We start from the inequality
$$
\int_t^T|f_s|\,dA_s = \int_t^Te^{-\frac{\beta}{2} A_s}e^{\frac{\beta}{2} A_s}|f_s|\,dA_s
\le
\left(  \int_t^T e^{-\beta A_s}  dA_s\right)^{1/2}\left(
  \int_t^T e^{\beta A_s}|{f}_s|^2 dA_s
\right)^{1/2}.
$$
Since $\beta\int_t^T e^{-\beta A_s}  dA_s= e^{-\beta A_t}-e^{-\beta A_T}\le e^{-\beta A_t}$ we
arrive at
\begin{equation}\label{felleunokappabeta}
\left(
\int_t^T|f_s|\,dA_s\right)^2 \le \frac{e^{-\beta A_t}}{\beta}
\int_t^T e^{\beta A_s}|{f}_s|^2 dA_s.
\end{equation}
This implies in particular that $\int_0^T|f_s|\,dA_s$ is square
summable. The solution $(Y,Z)$ is defined by considering a cadlag
version of the martingale $M_t=\E^{\calf_t}[\xi + \int_0^Tf_s  \,
dA_s ]$. By the martingale representation Theorem
\ref{rappresentazione}, there exists a  process $Z \in
\call^{1,0}(p)$ such that
\begin{equation*} M_t= M_0+ \int_0^t \int_K Z_s(y) \;q(dy \,ds) , \qquad t \in [0,  T].
 \end{equation*}
 Define the process $Y$ by
\begin{equation*}Y_t = M_t - \int_0^t f_s ( U_s,V_s) \; dA_s, \qquad t \in [0,  T].
 \end{equation*}
Noting that  $Y_T= \xi$, we easily deduce that the equation \eqref{BSDElineare} is satisfied.

It remains to show that $(Y,Z)\in\K^\beta$.
Taking the
conditional expectation, it follows from \eqref{BSDElineare} that
$Y_t= \E^{\calf_t}[\xi + \int_t^Tf_s  \, dA_s ]$ so that, using \eqref{felleunokappabeta}, we obtain
\begin{eqnarray}\nonumber
e^{\beta A_t} |Y_t|^2&\le& \dis
2e^{\beta A_t} |\E^{\calf_t}\xi|^2
+2e^{\beta A_t} \left|\E^{\calf_t}\int_t^Tf_s  \, dA_s\right|^2\\
&\le& \dis
 2\E^{\calf_t} \left[e^{\beta A_T}|\xi|^2 +  \frac{1}{\beta}\int_0^Te^{\beta A_s}|f_s|^2  \, dA_s\right].
 \label{martingausil}
 \end{eqnarray}
Denoting by $m_t$ the right-hand side of \eqref{martingausil}, we see that $m$ is a martingale
by the assumptions of the lemma. In particular, for every stopping time $S$ with values in $[0,T]$,
we have
\begin{equation}\label{uniftempiarresto}
\E\,e^{\beta A_S} |Y_S|^2\le \E\, m_S=\E\,m_T<\infty
\end{equation}
by the optional stopping theorem. Next we define the increasing sequence of stopping times
$$
S_n=\inf\{t\in [0,T]\,:\,
\int_0^t e^{\beta A_s }|Y_s|^2dA_s +
 \int_0^t  \int_Ke^{\beta A_s } |Z_s(y)|^2 \phi_s(dy)dA_s >n\},
 $$
with the  convention $\inf \emptyset =T$. Computing the Ito differential
$d( e^{\beta A_t}|{Y}_t|^2)$
on the interval $[0,S_n]$ and proceeding as before we deduce
$$
\begin{array}{l}\dis
\beta\, \E\int_0^{S_n}   e^{\beta A_s}| {Y}_s|^2 dA_s
+ \E\int_0^{S_n} \int_K e^{\beta A_s}| {Z}_s(y)|^2 \phi_s(dy) dA_s \\
\dis\qquad \le\E \,e^{\beta A_{S_n}}|Y_{S_n}|^2
  +2\E \int_0^{S_n}  e^{\beta A_s} {Y}_{s}  f_s\, dA_s.
\end{array}
$$
Using the inequalities
$2{Y}_{s}  f_s\le (\beta/2)|Y_s|^2 + (2/\beta)|f_s|^2$
and \eqref{uniftempiarresto} (with $S=S_n$) we find the following estimates
\begin{equation*}
   \E\int_0^{S_n}   e^{\beta A_s}| {Y}_s|^2 dA_s
 \le \frac{4}{\beta}\,
\E e^{\beta A_T}|\xi|^2 +  \frac{8}{\beta^2}\E\int_0^Te^{\beta A_s}|f_s|^2  \, dA_s,
\end{equation*}
\begin{equation*}
 \E\int_0^{S_n} \int_K e^{\beta A_s}| {Z}_s(y)|^2 \phi_s(dy) dA_s \le
4\E e^{\beta A_T}|\xi|^2 +  \frac{8}{\beta}\E\int_0^Te^{\beta A_s}|f_s|^2  \, dA_s,
\end{equation*}
from which we deduce
\begin{equation}\label{stimasonesseenne}
\begin{array}{l}\dis
   \E\int_0^{S_n}   e^{\beta A_s}| {Y}_s|^2 dA_s
+ \E\int_0^{S_n} \int_K e^{\beta A_s}| {Z}_s(y)|^2 \phi_s(dy) dA_s
\\\dis\qquad\qquad \le c_1(\beta)
\E e^{\beta A_T}|\xi|^2 +  c_2(\beta)\int_0^Te^{\beta A_s}|f_s|^2  \, dA_s,
\end{array}
\end{equation}
where $c_1(\beta)=4(1+\frac{1}{\beta})$ and $c_2(\beta)= \frac{8}{\beta}(1+\frac{1}{\beta})$.

Setting $S=\lim_nS_n$ we deduce
$$
\int_0^{S }   e^{\beta A_s}| {Y}_s|^2 dA_s
+  \int_0^{S } \int_K e^{\beta A_s}| {Z}_s(y)|^2 \phi_s(dy) dA_s <\infty,\qquad \P-a.s.
$$
which implies $S=T$, $\P$-a.s., by the definition of $S_n$. Letting $n\to\infty$
in \eqref{stimasonesseenne} we conclude that
\eqref{stimaprlineare} holds, so that $(Y,Z)\in \K^{\beta}$.
\qed

\begin{theorem} \label{Th:ex-un}
Suppose that Hypothesis \ref{hyp:BSDE} holds with
$\beta>L^2+2L'$.

Then there exists a unique pair $(Y,Z)$ in $\K^{\beta}$ which solves the BSDE \eqref{BSDE}.
\end{theorem}

\noindent{\bf Proof.}
We use a fixed point theorem for the mapping $\Gamma: \K^\beta \rightarrow \K^\beta$
 defined setting
 $(Y,Z)=\Gamma (U,V)$ if $(Y,Z)$ is the pair satisfying
\begin{equation}\label{mappabackward}
  Y_t+\int_t^T \int_K Z_s(y) \; q(ds \;dy)=
   \xi +\int_t^Tf_s ( U_s,V_s)\; dA_s.
\end{equation}
Let us remark that from the assumptions on $f$ it follows that $\E \int_0^T e^{\beta A_s}|f_s( U_s,V_s)|^2 dA_s < \infty$,
so by Lemma \ref{BSDEcasobanale} there exists a unique
$(Y,Z)\in\K^\beta$  satisfying
 \eqref{mappabackward} and
 $\Gamma$ is a
well defined map.

Let $(U^i,V^i)$, $i=1,2$, be elements of $\mathbb{K}^{\beta}$ and
let $(Y^i,Z^i)= \Gamma (U^i,V^i)$. Denote $\overline{Y}=Y^1-Y^2$,
$\overline{Z}=Z^1-Z^2$, $\overline{U}=U^1-U^2$,
$\overline{V}=V^1-V^2$, $\overline{f}_s=f_s(U^1_s,V^1_s)
-f_s(U^2_s,V^2_s)$. Lemma \ref{BSDEcasobanale} applies to
$\overline{Y},\overline{Z}, \overline{f}$ and
\eqref{identitaenergia} yields, noting that $\overline{Y}_T=0$,
$$
\begin{array}{l}\dis
\E e^{\beta A_t}| \overline{Y}_t|^2 +\beta\,\E\int_t^T   e^{\beta A_s}|\overline{Y}_s|^2 dA_s
+ \E\int_t^T \int_K e^{\beta A_s}|\overline{Z}_s(y)|^2 \phi_s(dy) dA_s \\
\dis\qquad
=
  2\E \int_t^T  e^{\beta A_s} \overline{Y}_{s} \,\overline{f}_s dA_s,
  \qquad t\in [0,T],
\end{array}
$$
From the Lipschitz conditions of $f$ and elementary inequalities it follows that
$$
\begin{array}{l}\dis
 \beta\,\E\int_0^T   e^{\beta A_s}|\overline{Y}_s|^2 dA_s
+ \E\int_0^T \int_K e^{\beta A_s}|\overline{Z}_s(y)|^2 \phi_s(dy) dA_s
\\
\dis\qquad\le
  2L\E \int_t^T  e^{\beta A_s} |\overline{Y}_{s}| \,
 \left( \int_K |\overline{V}_s(y)|^2 \phi_s(dy)\right)^{1/2}
   dA_s
  + 2L'\E \int_t^T  e^{\beta A_s} |\overline{Y}_{s}| \,|\overline{U}_s|\, dA_s
\\
\dis\qquad
\le
\alpha
\E \int_0^T \int_K e^{\beta A_s}|\overline{V}_s(y)|^2 \phi_s(dy)\, dA_s+
\frac{L^2}{\alpha} \E \int_0^T e^{\beta A_s}|\overline{Y}_s |^2 \, dA_s
\\
\dis\qquad
\quad+\gamma L' \E \int_0^T e^{\beta A_s}|\overline{Y}_s |^2 \, dA_s+
\frac{L'}{\gamma}  \E \int_0^T e^{\beta A_s}|\overline{U}_s |^2 \, dA_s
\end{array}
$$
for every $\alpha>0$, $\gamma>0$. This can be written
$$
\left( \beta- \frac{L^2}{\alpha} - \gamma L' \right)\,|\overline{Y}|^2_\beta
+ \|\overline{Z}\|_\beta^2 \le
\alpha
\|\overline{V}\|^2_\beta +
\frac{L'}{\gamma}  |\overline{U}  |^2_\beta.
$$
By the assumption on $\beta$ it is possible to find $\alpha\in (0,1)$ such that
$$
\beta>\frac{L^2}{\alpha} +\frac{2L'}{\sqrt\alpha}.
$$
If $L'=0$ we see that $\Gamma$   is an $\alpha$-contraction on
$\K^\beta$ endowed with the equivalent norm
$(Y,Z)\mapsto ( \beta- ({L^2}/{\alpha}))\,|{Y}|^2_\beta
+ \|{Z}\|_\beta^2$. If $L'>0$ we choose $\gamma=1/\sqrt{\alpha}$
and obtain
$$
\frac{L'}{\sqrt\alpha} \,|\overline{Y}|^2_\beta
+ \|\overline{Z}\|_\beta^2 \le
\alpha
\|\overline{V}\|^2_\beta +
 {L'}{\sqrt{\alpha}}  |\overline{U}  |^2_\beta=
 \alpha \,
 \left( \frac{L'}{\sqrt\alpha} \,|\overline{U}|^2_\beta
+ \|\overline{V}\|_\beta^2 \right),
$$
so that
$\Gamma$   is an $\alpha$-contraction on
$\K^\beta$ endowed with the equivalent norm
$(Y,Z)\mapsto (  {L'}/\sqrt{\alpha} )\,|{Y}|^2_\beta
+ \|{Z}\|_\beta^2$. In all cases
there exists a unique fixed point which is the required unique solution to the BSDE \eqref{BSDE}.
\qed

We next prove some estimates on the solutions of the BSDE, which
show in particular the continuous dependence upon the data. Let us
consider two solutions $(Y^1,Z^1)$, $(Y^2,Z^2)\in \K^{\beta}$ to
the BSDE \eqref{BSDE} associated with the drivers $f^1$ and $f^2$
and final data $\xi^1$ and $\xi^2$, respectively, which are
assumed to satisfy Hypothesis \ref{hyp:BSDE}. Denote
$\overline{Y}=Y^1-Y^2$, $\overline{Z}=Z^1-Z^2$, $\overline{\xi}=
\xi^1-\xi^2$, $ \overline{f}_s= f^1_s(Y^2_s,Z^2_s(\cdot))-
f^2_s(Y^2_s,Z^2_s(\cdot)).
 $


\begin{proposition} \label{a priori estimates}
Let $(\overline{Y},\overline{Z})$ be the processes defined above.
Then, for $\beta > 2L' + L^2$, the a priori estimates hold:
\begin{eqnarray*}\label{stima-diff-y}
|\overline{Y}|^2_{\beta} & \leq & \frac{2}{\beta-2L' -L^2}\,\E e^{\beta A_T}|\overline{\xi}|^2 + \frac{4}{(\beta-2L' -L^2)^2}\,
\E \int_0^T e^{\beta A_s}|\overline{f}_s|^2 dA_s,\\
\label{stima-diff-z}
\|\overline{Z}\|^2_{\beta} & \leq  & \left(2+ \frac{16}{\beta-2L' -L^2}\right)\,\E e^{\beta A_T}|\overline{\xi}|^2
\\
&&+ \frac{2}{\beta-2L' -L^2}\left(1+\frac{16}{\beta-2L' -L^2}\right)\,\E \int_0^T e^{\beta A_s}|\overline{f}_s|^2 dA_s.
\end{eqnarray*}
\end{proposition}

\noindent{\bf Proof.}
From the Ito formula applied to $ e^{\beta A_t}|\overline{Y}_t|^2$ it follows that
$$d( e^{\beta A_t}|\overline{Y}_t|^2)= \beta  e^{\beta A_t}|\overline{Y}_t|^2 dA_t + 2 e^{\beta A_t}\overline{Y}_{t-} dY_t + e^{\beta A_t}|\Delta \overline{Y}_t|^2.$$

So integrating on $[t,T]$ and recalling that $A$ is continuous,
\begin{eqnarray}
\nonumber e^{\beta A_t}|\overline{Y}_t|^2 & = & - \int_t^T \beta  e^{\beta A_s}|\overline{Y}_s|^2 dA_s - 2 \int_t^T e^{\beta A_s}\overline{Y}_{s-} \int_K \overline{Z}_s(y) q(ds\,dy) - \sum_{t<s\le T} e^{\beta A_s}|\Delta \overline{Y}_s|^2\\
\label{eq:int-Ito} &  & +e^{\beta A_T}|\overline{\xi}|^2
+2\int_t^T  e^{\beta A_s}\overline{Y}_{s} (f^1(Y_s^1,Z^1_s(\cdot))-f^2(Y_s^2,Z^2_s(\cdot)) dA_s.
\end{eqnarray}
The integral process $\int_0^t  e^{\beta A_s}\overline{Y}_{s-} \int_K Z_s(y) q(ds\,dy)$
is a martingale, because the integrand process $ e^{\beta A_s}\overline{Y}_{s-}  \overline{Z}_s(y)$ is in $L^1(p)$:
 in fact from the Young inequality we get
\begin{eqnarray*}
             \E \int_0^T \int_K  e^{\beta A_s}|\overline{Y}_{s-}| | \overline{Z}_s (y)| \phi_s(dy)dA_s && \\
     \le \frac{1}{2}\E \int_0^T e^{\beta A_s}|\overline{Y}_{s-}|^2 dA_s &+&
      \frac{1}{2}\E \int_0^T \int_K  e^{\beta A_s}| \overline{Z}_s (y)|^2 \phi_s(dy)dA_s < +\infty.
             \end{eqnarray*}

Moreover  we have
\begin{eqnarray*}
\sum_{0<s\le t} e^{\beta A_s}|\Delta \overline{Y}_s|^2 & = & \int_0^t \int_K e^{\beta A_s}|\overline{Z}_s(y)|^2p(ds\,dy)
  \\
 & = &\int_0^t \int_K e^{\beta A_s}|\overline{Z}_s(y)|^2 q(ds\,dy)+
  \int_0^t \int_K e^{\beta A_s}|\overline{Z}_s(y)|^2 \phi_s(dy) dA_s,
 \end{eqnarray*}
where the stochastic integral with respect to $q$ is a martingale.

Hence taking the expectation in \eqref{eq:int-Ito}, by the
Lipschitz property of the driver $f^1$ and using the notation
$\|z(\cdot)\|_s^2=\int_K |z(y)|^2\phi_s(dy)$ we get
\begin{eqnarray*}
\E e^{\beta A_t}|\overline{Y}_t|^2 &=& -\E\int_t^T \beta  e^{\beta A_s}|\overline{Y}_s|^2 dA_s - \E\int_t^T \int_K e^{\beta A_s}|\overline{Z}_s(y)|^2 \phi_s(dy) dA_s +\E e^{\beta A_T}|\overline{\xi}|^2\\
& &+2\E \int_t^T  e^{\beta A_s}\overline{Y}_{s} (f^1(Y_s^1,Z^1_s)-f^2(Y_s^2,Z^2_s) dA_s\\
& \leq & -\E\int_t^T \beta  e^{\beta A_s}|\overline{Y}_s|^2 dA_s
- \E\int_t^T \int_K e^{\beta A_s}|\overline{Z}_s(y)|^2 \phi_s(dy) dA_s +\E e^{\beta A_T}|\overline{\xi}|^2\\
&   &+2\E \int_t^T  e^{\beta A_s}|\overline{Y}_{s}|(|f^1(Y_s^1,Z^1_s)-f^1(Y_s^2,Z^2_s)| +|\overline{f}_s|)dA_s\\
& \leq & -\E\int_t^T \beta  e^{\beta A_s}|\overline{Y}_s|^2 dA_s
- \E\int_t^T e^{\beta A_s}\|\overline{Z}_s \|_s^2   dA_s
+\E e^{\beta A_T}|\overline{\xi}|^2\\
&  &+2L'\E \int_t^T \! e^{\beta A_s}|\overline{Y}_{s}|^2 dA_s\!
+\! 2L\E \int_t^T  \!e^{\beta
A_s}|\overline{Y}_{s}|\|\overline{Z}_s\|_sdA_s \!+\!2\E \int_t^T
\!e^{\beta A_s}|\overline{Y}_{s}||\overline{f}_s|dA_s.
\end{eqnarray*}
We note that the quantity $Q(y,z)=-\beta|y|^2 - \|z\|_s^2+ 2L'|y|^2 +2L|y|\|z\|_s +2|\overline{f}_s||y| $
which occurs in the integrand terms in the right hand of the above inequality can be written as
\begin{eqnarray*}
Q(y,z) & = & -\beta|y|^2+ 2L'|y|^2 +L^2|y|^2  +2|\overline{f}_s||y| - (\|z\|_s- L|y|)^2\\
 & = & -\beta_L(|y|- \beta_L^{-1}|\overline{f}_s|)^2 - (\|z\|_s- L|y|)^2 +\beta_L^{-1}|\overline{f}_s|^2
\end{eqnarray*}
where $\beta_L:= \beta -2L'-L^2$ is assumed to be strictly positive.
Hence
\begin{eqnarray*}\E e^{\beta A_t}|\overline{Y}_t|^2  &+& \beta_L \E\int_t^T e^{\beta A_s}(|\overline{Y}_s|- \beta_L^{-1}|\overline{f}_s|)^2 dA_s+ \E\int_t^T  e^{\beta A_s}(\|\overline{Z}_s\|_s -L|\overline{Y}_s|)^2  dA_s \\
&\leq& \E e^{\beta A_T}|\overline{\xi}|^2 +\E \int_t^T  e^{\beta A_s}\frac{|\overline{f}_s|^2}{\beta_L} dA_s
\end{eqnarray*}
from which we deduce
\begin{eqnarray*}\E \int_0^T e^{\beta A_s}|\overline{Y}_s|^2 dA_s &\leq& \frac{2}{\beta_L}\E e^{\beta A_T}|\overline{\xi}|^2 +\frac{4}{\beta_L^2}\E \int_0^T  e^{\beta A_s}|\overline{f}_s|^2 dA_s
\end{eqnarray*}
and
\begin{eqnarray*} \E\int_0^T  e^{\beta A_s}\|\overline{Z}_s\|_s ^2  dA_s &\leq& \left(2 +\frac{16}{\beta_L}\right)\E e^{\beta A_T}|\overline{\xi}|^2 +\frac{2}{\beta_L}\left(1+\frac{16}{\beta_L}\right)\E \int_0^T  e^{\beta A_s}\frac{|\overline{f}_s|^2}{\beta_L} dA_s
\end{eqnarray*}
 \qed

From the a priori estimates one can deduce the continuous dependence of the solution upon the data.
\begin{proposition}
Suppose that Hypothesis \ref{hyp:BSDE} holds with
$\beta>L^2+2L'$ and
let $(Y,Z)$ be the unique solution in $\K^{\beta}$ to the BSDE (\ref{BSDE}).
Then
\begin{eqnarray*}
\E\int_0^T e^{\beta A_s}|Y_s|^2 dA_s  &+&  \E\int_0^T  e^{\beta A_s}\int_K |Z_s(y)|^2\phi_s(dy)  dA_s
\\
 &\leq &C_1(\beta)\E e^{\beta A_T}|\xi|^2 +C_2(\beta)\int_0^T  e^{\beta A_s}\frac{|\overline{f}_s|^2}{\beta_L} dA_s,
\end{eqnarray*}
where
$C_1(\beta)=\left(2 +\frac{18}{\beta -2L'-L^2}\right)$, $C_2(\beta)=\frac{2}{\beta -2L'-L^2}\left(1+\frac{18}{\beta -2L'-L^2}\right)$.
\end{proposition}

\noindent{\bf Proof.}
The thesis follows from Proposition \ref{a priori estimates} setting $f^1=f$, $\xi^1=\xi$, $f^2=0$ and $\xi^2=0$.
\qed

\section{Optimal control}\label{sec-control}

Throughout this section we assume that
a marked point process is given, satisfying the assumptions of Section \ref{notations}.
In particular we suppose that $T_n\to\infty$ $\P$-a.s. and that
(\ref{acontinuo}) holds.

The data specifying the optimal control problem are an
  action (or decision) space
$U$, a running cost function $l$, a terminal cost function $g$,
and another function $r$ specifying the effect of the control
process. They are assumed to satisfy the following conditions.

\begin{hypothesis}\label{hyp:controllo}
\begin{enumerate}
\item  $(U,\calu)$ is a measurable space.
\item  The functions $r,l:\Omega\times [0,T]\times K\times U\to \R$
are $\calp\otimes \calk\otimes \calu$-measurable and
 there exist
constants $C_r> 1$, $C_l>0$ such that, $\P$-a.s.,
\begin{equation}\label{ellelimitato}
0\le r_t (y,u)\le C_r,\quad |l_t (x,u)|\le C_l,\qquad t\in [0,T],
x\in K, u\in U.
\end{equation}
\item
The function $g:\Omega\times K\to \R$ is $\calf_T\otimes
\calk$-measurable.
\end{enumerate}
\end{hypothesis}

We define as an admissible  control process, or simply a control,
any predictable process $(u_t)_{t\in[0,T]}$ with values in $U$. The set of admissible
control processes is denoted $\cala$.

To every control $u(\cdot)\in\cala$ we associate a probability measure $\P_u$ on
$(\Omega,\calf)$ by a change of measure of Girsanov type, as we now
describe. We define
$$
L_t=
\exp\left(\int_0^t\int_K (1-r_s (y,u_s))\;\phi_s(dy)\,dA_s\right)
\prod_{n\ge1\,:\,T_n\le t}r_{T_n} (\xi_n,u_{T_n}),
\qquad t\in [0,T],
$$
with the convention that the last product equals $1$ if there are
no indices $n\ge 1$ satisfying $T_n\le t$ (similar conventions
will be adopted later without further mention). It is a well-known
result that $L$ is a nonnegative supermartingale, (see \cite{J}
Proposition 4.3, or \cite{Bo-Va-Wo}), solution to the equation
$$
L_t=
1+\int_0^t\int_K L_{s-}\,(r_s (y,u_s)-1)\;q(ds\,dy),
\qquad t\in [0,T].
$$
The following result  collects some properties
of the process $L$ that we need later.

\begin{lemma}\label{gammasommabile}
Let $\gamma>1$ and
\begin{equation}\label{valoredibeta}
    \beta= \gamma+1+\frac{C_r^{\gamma^2}}{\gamma-1}.
\end{equation}
If $\E \exp (\beta A_T)<\infty$ then  we have $\sup_{t\in [0,T]}\E L_t^\gamma<\infty$ and
$\E L_T=1$.
\end{lemma}

\noindent{\bf Proof.} We follow \cite{B}, Chapter VIII Theorem T11, with some
modifications.
To shorten notation we define $\rho_s(y)=r_s (y,u_s)$ and we denote $L_t=\cale (\rho)_t$. For
$\gamma>1$ we define
$$
a_s(y)=\gamma^{-1}(1-\rho_s(y)^{\gamma^2}),
\qquad
b_s(y)=\gamma-\gamma\rho_s(y) -\gamma^{-1}+ \gamma^{-1} \rho_s(y)^{\gamma^2},
$$
so that $\gamma (1-\rho_s (y))=a_s(y)+b_s(y)$.
Then
$$
L_t^\gamma=
\exp\left(\int_0^t\int_K (a_s(y)+b_s(y))\,\phi_s(dy)\,dA_s\right)
\prod_{T_n\le t}\rho_{T_n} (\xi_n)^\gamma,
$$
and by H\"older's inequality
$$\begin{array}{lll}
\E L_t^\gamma&\le &\dis \left\{\E\bigg[
\exp\left(\int_0^t\int_K \gamma a_s(y)\,\phi_s(dy)\,dA_s\right)
\prod_{T_n\le t}\rho_{T_n} (\xi_n)^{\gamma^2}\bigg]
\right\}^{\frac1\gamma}
\\&&\dis
\left\{\E
\exp\left(\int_0^t\int_K \frac{\gamma}{\gamma-1} b_s(y) \;\phi_s(dy)\,dA_s\right)
\right\}^{\frac{\gamma-1}{\gamma}}.
\end{array}
$$
Noting that  $\gamma a_s(y)=1-\rho_s(y)^{\gamma^2}$, the  term in square brackets
equals $\cale (\rho^{\gamma^2})_t$ and we have $\E \cale (\rho^{\gamma^2})_t\le 1$ by
the supermartingale property. Since
$
b_s(y)\le \gamma  -\gamma^{-1}+ \gamma^{-1} C_r^{\gamma^2}$ we arrive at
\begin{equation}\label{ellegammaintegrabile}
    \E L_t^\gamma \le
\left\{\E
\exp\left(A_T \bigg(\gamma+1+\frac{C_r^{\gamma^2}}{\gamma-1} \bigg)\right)
\right\}^{\frac{\gamma-1}{\gamma}}
=
\left\{\E
\exp\left(\beta A_T\right)
\right\}^{\frac{\gamma-1}{\gamma}}<\infty
.
\end{equation}

Let $S_n=\inf \{t\in [0,T]\,:\, L_{t-}+A_t\ge n\}$ with the convention $\inf \emptyset =T$,
and let $\rho_s^{(n)}(y)=1_{[0,S_n]}(s)\rho_s(y)+1_{(S_n,T]}(s)$, $L^{(n)}=\cale (\rho^{(n)})$.
Then $L^{(n)}$ satisfies
$$
L_t^{(n)}=
1+\int_0^t\int_K L_{s-}^{(n)}\,(r_s^{(n)} (y)-1)\;q(ds\,dy),
\qquad t\in [0,T].
$$
By the choice of $\rho^{(n)}$ we have $L^{(n)}_t=L_{t\wedge S_n}$,
and by the choice of $S_n$ it is easily proved that
$\E\int_0^T\int_K L_{s-}^{(n)}\,|r_s^{(n)} (y)-1|\,\phi_s(dy)\,dA_s<\infty$, so that
$L^{(n)}$ is a martingale and $\E L^{(n)}_t=\E L_{t\wedge S_n}=1$.
The first part of the proof applies to $L^{(n)}$ and
the inequality (\ref{ellegammaintegrabile})
yields in particular
$\sup_n \E (L^{(n)}_{t})^\gamma=\sup_n \E (L_{t\wedge S_n})^\gamma <\infty$.
So $(L_{t\wedge S_n})_n$ is uniformly integrable and
letting $n\to\infty$ we conclude that
$\E L_t=1$.
\qed

Under the assumptions of the lemma, the process $L$ is a
martingale and we can define a probability $\P_u$ setting
$\P_u(d\omega)=L_T(\omega)\P(d\omega)$. It can then be proved (see
\cite{J} Theorem 4.5) that the compensator $\tilde{p}^u$ of $p$
under $\P_u$ is related to the compensator $\tilde{p}$ of $p$
under $\P$ by the formula
$$
\tilde{p}^u(dt\,dy)=r_t(y,u_t)\,\tilde{p}(dt\,dy)=r_t(y,u_t)\,\phi_t(dy)\,dA_t.
$$
In particular the compensator of $N$ under $\P_u$ is
\begin{equation}\label{aucompensato}
    A^u_t=\int_0^t\int_K r_s(y,u_s)\,\phi_s(dy)\,dA_s.
\end{equation}

We finally define the cost  associated to every $u(\cdot)\in\cala$
as
$$
J(u(\cdot))=\E_u\left[
\int_0^Tl_t(X_t,u_t)\;dA_t + g(X_T)\right],
$$
where $\E_u$ denotes the expectation under $\P_u$. Later we will
assume that
\begin{equation}
\label{condsucostofinale}
 \E[|g(X_T)|^2e^{\beta A_T} ]<\infty
 \end{equation}
  for some $\beta>0$
that will be fixed in such a way that the cost is finite for every
admissible control. The control problem consists in minimizing $
J(u(\cdot))$ over $\cala$.

\begin{remark}\begin{em}
We recall (see e.g. \cite{B}, Appendix A2, Theorem T34) that a
process $u$ is $(\calf_t)$-predictable if and only if it admits
the representation
\begin{equation}\label{rappresentazione-predictable-proc}
u(\omega,t)= \sum_{n \geq 0} u^{(n)}(\omega,t)\,1_{T_n(\omega) <t
\leq T_{n+1}(\omega)}
\end{equation}
where for each $n \ge 0$ the mapping $(\omega,t) \mapsto
u^{(n)}(\omega,t)$ is $\calf_{T_n} \otimes \calb
([0,\infty))$-measurable. Since we have
$\calf_{T_n}=\sigma(T_i,\xi_i, 0\le i\le n)$ (see e.g. \cite{B},
Appendix A2, Theorem T30) the fact that a control is predictable
can be roughly interpreted by saying that the controller, at each
time  $T_n$, based on observation of the random variables
$T_i,\xi_i, 0\le i\le n$,  chooses his present and future control
actions and updates his decisions only at time $T_{n+1}$.
\end{em}
\end{remark}

\begin{remark}\begin{em}
We notice that the laws of the random coefficients $r,l,g$ under
$\P$ and under $\P_u$  are not the same in general, so that the
formulation of the optimal control problem should be carefully
examined when facing a specific application or modeling situation.
This difficulty clearly disappears when $r,l,g$ are deterministic.
\end{em}
\end{remark}

We next proceed to the solution of the optimal control problem
formulated above. A basic role is played by  the BSDE
\begin{equation}\label{bsdecontrollo}
    Y_t+\int_t^T\int_KZ_s(y)\,q(ds\,dy) =
g(X_T) +\int_t^T f(s,X_s,Z_s(\cdot))\,dA_s, \qquad t\in [0,T],
\end{equation}
with terminal condition $g(X_T)$ and generator defined by means of
the hamiltonian function $f$. The hamiltonian function is defined
for every $\omega\in\Omega, t\in[0,T], x\in K$ and $z\in
\call^1(K,\calk, \phi_t(\omega,dy))$ by the formula
\begin{equation}\label{defhamiltonian}
    f(\omega,t,x,z(\cdot))=\inf_{u\in U}\left\{
l_t(\omega, x,u)+ \int_K z(y) \, (r_t
(\omega,y,u)-1)\,\phi_t(\omega,dy)\right\}.
\end{equation}

We will assume that the infimum is in fact achieved, possibily at
many points. Moreover we need to verify that the generator of the
BSDE satisfies the conditions required in the previous section. It
turns out that an appropriate assumption is the following one, since we
will see below (compare Proposition \ref{selezione}) that it can
be verified under quite general conditions. Here and in the following we set $X_{0-}=X_0$.

\begin{hypothesis}\label{hyp:hamiltoniana}
For every $Z\in\call^{1,0}(p)$ there exists a function
$\underline{u}^Z:\Omega\times [0,T]\to U$, measurable with
respect to $\calp$ and $\calu$, such that
\begin{equation}\label{minselector}
    f(\omega,t,X_{t-}(\omega),Z_t(\omega,\cdot))=
l_t(X_{t-}(\omega),\underline{u}^Z(\omega,t )) + \int_K Z_t(\omega, y) \,
\Big(r_t
(\omega,y,\underline{u}^Z(\omega,t))-1\Big)\,\phi_t(\omega,dy)
\end{equation}
for almost
all $(\omega,t)$ with respect to the measure
$dA_t(\omega)\P(d\omega)$.
\end{hypothesis}
 Note that if $ Z \in
\mathcal{L}^{1,0}(p)$ then $ Z_t(\omega,\cdot)$ belongs to
$\call^1(K,\mathcal{K}, \phi_t(\omega,dy) )$ except possibly on a
predictable set of points $(\omega,t )$ of measure zero with
respect to $dA_t(\omega)\P(d\omega)$, so that the equality
\eqref{minselector} is meaningful.
Also note that each $u^Z$ is an admissible control.

We can now verify that all the assumptions of Hypothesis
\ref{hyp:BSDE} hold true for the generator of the BSDE
\eqref{bsdecontrollo}, which is given by the formula
$$
f_t(\omega, z(\cdot))= f(\omega,t,X_t(\omega),z(\cdot)), \qquad
\omega \in \Omega, \;t \in [0,T],\; z\in\call^{2}(K,
\calk,\phi_t(\omega,dy)).
$$
Indeed, if $ Z\in\mathcal{L}^{2,\beta}(p)$ then $Z\in
\mathcal{L}^{1,0}(p)$ by \eqref{inclusioneprocessi}, and
\eqref{minselector} shows that the process $(\omega,t)\mapsto
f(\omega,t,X_{t-}(\omega),Z_t(\omega,\cdot))$ is progressive;
since $A$ is assumed to have continuous trajectories and $X$ has piecewise
constant pahts, the progressive set
$\{(\omega,t)\,:\, X_{t-}(\omega)\neq X_{t}(\omega)\}$ has
measure zero with respect to $dA_t(\omega)\P(d\omega)$;  it follows that the
process
$$(\omega,t)\mapsto f(\omega,t,X_t(\omega),Z_t(\omega,\cdot))
=f_t(\omega, Z_t(\omega,\cdot))
$$
is progressive, after modification on a set of measure zero,
as required in \eqref{misurabilitaf}. Next, using
the boundedness assumptions \eqref{ellelimitato},
 it is easy to check  that \eqref{generatorelip} is verified with
$L'=0$ and
$$L=\sup|r-1|=\sup\{ |r_t (y,u)-1|\,:\,
\omega\in\Omega,\,  t\in [0,T], \,y\in  K, \,u\in U\}.
$$
Using \eqref{ellelimitato} again we also have
\begin{equation}\label{gensommcontrollo}
  \E \int_0^T e^{\beta A_t}|f(t,X_t,0)|^2 dA_t
  = \E \int_0^T e^{\beta A_t}|\inf_{u\in U}  l_t( X_t,u)|^2 dA_t
  \le C_l^2\,\beta^{-1}\,
\E \,e^{\beta A_T},
\end{equation}
 so  that \eqref{generatoresommabile} holds as well provided the
 right-hand side of \eqref{gensommcontrollo}
 is finite.
 Assuming finally that \eqref{condsucostofinale} holds, by
Theorem \ref{Th:ex-un} the BSDE has a unique solution
$(Y,Z)\in\K^\beta$ if $\beta>L^2$.

The corresponding admissible control $\underline{u}^Z$, whose
existence is required in Hypothesis \ref{hyp:hamiltoniana}, will be denoted
$ u^*$.

 We are now ready to state the main result of this section.
Recall that $C_r> 1$ was introduced in \eqref{ellelimitato}.

\begin{theorem}\label{teoremacontrollononmarkov}
Assume that Hypotheses \ref{hyp:controllo} and
\ref{hyp:hamiltoniana} are satisfied
  and that
\begin{equation}\label{ATintegrabile}
\E \exp \left((3+C_r^4) A_T\right)<\infty.
\end{equation}
Suppose also that there exists $\beta$ such that
\begin{equation}\label{compatibbeta}
  \beta>\sup|r-1|^2, \quad \E \exp \left(\beta A_T\right)<\infty,
\quad \E[|g(X_T)|^2e^{\beta A_T} ]<\infty.
\end{equation}
Let $(Y,Z)\in \K^\beta$ denote
 the solution to the BSDE \eqref{bsdecontrollo} and
$u^*=\underline{u}^Z$ the corresponding admissible control.
Then
$u^*(\cdot)$ is   optimal and $Y_0$ is the optimal cost, i.e.
$Y_0= J(u^*(\cdot))= \inf_{u(\cdot)\in\cala }J(u(\cdot))$.
\end{theorem}

\begin{remark} \begin{em}
Note that if $g$ is bounded then \eqref{compatibbeta} follows from
\eqref{ATintegrabile} with
$\beta=3+C_r^4$, since $|r_t (y,u)-1|^2\le (C_r+1)^2<3+C_r^4$.
\end{em}
\end{remark}

 \noindent{\bf Proof.} Fix $u(\cdot)\in\cala$.
Assumption \eqref{ATintegrabile} allows to apply Lemma
\ref{gammasommabile} with
 $\gamma=2$ and yields $\E L_T^2<\infty$.
 It follows that $g(X_T)$ is integrable under $\P_u$. Indeed by \eqref{condsucostofinale}
 $$
\E_u |g(X_T)|= \E |L_Tg(X_T)|\le (\E L_T^2)^{1/2}(\E g(X_T)^2)^{1/2}<\infty.
$$

 We next  show that under $\P_u$ we have $Z\in \call^{1,0}(p)$, i.e.
 $\E_u \int_{0}^T \int_K |Z_t(y)|\,\tilde{p}^u(dt\,dy)
<\infty$.
First note that, by H\"older's inequality,
$$
\begin{array}{lll}\dis
\int_{0}^T \int_K  |Z_t(y)|\,\phi_t(dy)\,dA_t&=&\dis
  \int_{0}^T \int_K e^{-\frac{\beta}{2}A_t}e^{\frac{\beta}{2}A_t} |Z_t(y)|\,\phi_t(dy)\,dA_t
\\
&\le&\dis
\left( \int_{0}^T   e^{-\beta A_t} dA_t\right)^{1/2}
\left(
\int_{0}^T \int_K e^{\beta A_t} |Z_t(y)|^2\,\phi_t(dy)\,dA_t
 \right)^{1/2}
\\
&=&\dis
\left( \frac{1-   e^{-\beta A_T}}{\beta} \right)^{1/2}
\left(
\int_{0}^T \int_K e^{\beta A_t} |Z_t(y)|^2\,\phi_t(dy)\,dA_t
 \right)^{1/2}.
\end{array}
$$
Therefore, using \eqref{ellelimitato},
$$
\begin{array}{lll}\dis
\E_u \int_{0}^T \int_K |Z_t(y)|\,\tilde{p}^u(dt\,dy)&=&\dis
\E_u \int_{0}^T \int_K  |Z_t(y)|\, r_t(y,u_t) \,\phi_t(dy)\, dA_t
\\&
=&\dis
\E\left[L_T \int_{0}^T \int_K  |Z_t(y)|\, r_t(y,u_t) \,\phi_t(dy)\, dA_t\right]
\\&
\le&\dis
(\E L_T^2)^{1/2}\frac{C_r}{\sqrt{\beta}}
\left\{\E\int_{0}^T \int_K  e^{\beta A_t}|Z_t(y)|^2\,\phi_t(dy)\,dA_t\right\}^{1/2}
\end{array}
$$
and the right-hand side of the last inequality is finite,
since $(Y,Z)\in\K^\beta$. We have now proved that
 $Z\in \call^{1,0}(p)$ under $\P_u$.

 In particular it follows that
$$
\E_u \int_{0}^T \int_K Z_t(y)\,{p}(dt\,dy)=
\E_u \int_{0}^T \int_K Z_t(y)\,\tilde{p}^u(dt\,dy)=
\E_u \int_{0}^T \int_K   Z_t(y) \, r_t(y,u_t) \,\phi_t(dy)\, dA_t.
$$
Setting $t=0$ and taking the expectation $\E_u$
in the BSDE (\ref{bsdecontrollo}), recalling that
$q(dt\,dy)=p(dt\,dy)-\tilde{p}(dt\,dy)= p(dt\,dy)-\phi_t(dy)\,dA_t$
and that $Y_0$
is deterministic,
we obtain
$$
    Y_0+
    \E_u \int_{0}^T \int_K   Z_t(y) \, (r_t(y,u_t)-1) \,\phi_t(dy)\, dA_t=
\E_u \,g(X_T) +\E_u\int_0^T f(t,X_t,Z_t(\cdot))\,dA_t.
$$
We finally obtain
\begin{eqnarray*}
      Y_0 &=& J(u(\cdot))+
    \E_u \int_{0}^T\left[ f(t,X_t,Z_t(\cdot))
    -l_t(X_t,u_t)
    -\int_K   Z_t(y) \, (r_t(y,u_t)-1) \,\phi_t(dy)\,
    \right]\,dA_t \\
    &=& J(u(\cdot))+
    \E_u \int_{0}^T\left[ f(t,X_{t-},Z_t(\cdot))
    -l_t(X_{t-},u_t)
    -\int_K   Z_t(y) \, (r_t(y,u_t)-1) \,\phi_t(dy)\,
    \right]\,dA_t,
\end{eqnarray*}
where the last equality follows from the continuity if $A$.
This identity is sometimes called the fundamental relation. By the definition
of the hamiltonian $f$, the term in square brackets is smaller or equal to $0$,
and it equals $0$ if $u(\cdot)=u^*(\cdot)$.
\qed

Hypothesis \ref{hyp:hamiltoniana} can be verified in specific
situations when it is possible to compute explicitly the functions
$\underline{u}^Z$. General conditions for its validity can also be
formulated using appropriate selection theorems, as in  the
following proposition.

\begin{proposition} \label{selezione}
In addition to the assumptions in Hypothesis \ref{hyp:controllo},
suppose that $U$ is a compact metric space  with its
Borel $\sigma$-algebra $\calu$ and that the functions
$r_t(\omega,x,\cdot), l_t(\omega,x,\cdot):U\to \R$ are
continuous for every $\omega\in\Omega$, $t\in [0,T]$, $x\in
K$. Then Hypothesis \ref{hyp:hamiltoniana} is verified.
\end{proposition}

\noindent {\bf Proof.}
Let us consider the measure
$\mu(d\omega\,dt)=dA_t(\omega)\P(d\omega)$ on the predictable
$\sigma$-algebra $\calp$.  Let $\bar\calp$ denote its $\mu$-completion and
consider the complete measure space $(\Omega\times [0,T],  \bar\calp, \mu)$.
 Fix $Z\in\call^{1,0}(p)$, note that the
set $A^Z=\{(\omega,t)\,:\, Z_t(\omega,
\cdot)\notin\call^1(K,\calk,\phi_t(\omega,dy))$ has $\mu$-measure zero
 and define a map $F^Z:\Omega\times [0,T] \times U\to\R $ setting
$$
F^Z(\omega, t,u)=\left\{\begin{array}{ll} \dis l_t(\omega, X_{t-}(\omega),u)
+ \int_K Z_t(\omega, y) \, \Big(r_t
(\omega,y,u)-1\Big)\,\phi_t(\omega,dy) & \text{ if } (\omega,
t)\notin A^Z,
\\
0   & \text{ if } (\omega, t)\in A^Z.
\end{array}
\right.
$$
Then $ F^Z(\cdot,\cdot, u)$ is $\bar\calp$-measurable for every $u\in U$, and
it is easily verified that $ F^Z(\omega, t,\cdot)$
is continuous for every $(\omega,t)\in \Omega\times [0,T]$.
By a
classical selection theorem (see \cite{AuFr}, Theorems 8.1.3 and 8.2.11)
 there exists a function
$\underline{u}^Z:\Omega\times [0,T]\to U$, measurable with respect to $\bar\calp$ and
$\calu$, such that $F^Z(\omega,t,\underline{u}^Z(\omega,t)) =
\min_{u\in U} F^Z(\omega,t,u)$ for every $(\omega,t)\in \Omega\times [0,T]$, so that
 \eqref{minselector} holds true for every $(\omega,t)$.  After modification on a set of $\mu$-measure
zero, the function $u^Z$ can be made
measurable with respect to $\calp$ and
$\calu$, and \eqref{minselector} still holds, as it is understood as an equality
for $\mu$-almost
all $(\omega,t)$.
\qed

In several contexts, for instance in order to apply dynamic
programming arguments, it is useful to introduce a family of
control problems parametrized by $(t,x)\in [0,T]\times K$. Recall
the definition of the processes $(X_s^{t,x})_{s\in [t,T]}$ in
subsection \ref{famigliaparametrizzata}.

For fixed $(t,x)$ the cost corresponding to $u\in\cala$ is
defined as the random variable
$$
J_t(x,u(\cdot))=\E_u^{\calf_t}\left[
\int_t^Tl_s(X_s^{t,x},u_s)\;dA_s + g(X_T^{t,x})\right],
$$
where $\E_u^{\calf_t}$ denotes the conditional
expectation under $\P_u$ given ${\calf_t}$. We also introduce the (random)
value function
$$
v(t,x)= \mathop{\rm ess\ inf}_{u(\cdot)\in\cala }J_t(x,u(\cdot)),
\qquad t\in [0,T], \,x\in K.
$$

For every $(t,x)\in [0,T]\times K$ we
 consider the BSDE
\begin{equation}\label{bsdecontrollobis}
    Y_s^{t,x}+\int_s^T\int_KZ_r^{t,x}(y)\,q(dr\,dy) =
g(X_T^{t,x}) +\int_s^T f(r,X_r^{t,x},Z_r^{t,x}(\cdot))\,dA_r,
\qquad s\in [t,T].
\end{equation}
We need the following extended variant of Hypothesis \ref{hyp:hamiltoniana},
where we set $X_{t-}^{t,x}=x$:

\begin{hypothesis}\label{hyp:hamiltonianabis} For every $(t,x)\in [0,T]\times K$
and every $Z\in\call^{1,0}(p)$ there exists a function
$\underline{u}^{Z,t,x}:\Omega\times [t,T]\to U$, measurable with
respect to $\calp$ and $\calu$, such that
$$
    f(\omega,s,X_{s-}^{t,x}(\omega),Z_s(\omega,\cdot))=
l_t(X_{s-}(\omega),\underline{u}^{Z,t,x}(\omega,s )) + \int_K Z_s(\omega, y) \,
\Big(r_s
(\omega,y,\underline{u}^{Z,t,x}(\omega,s))-1\Big)\,\phi_s(\omega,dy)
$$
for almost
all $(\omega,s)\in \Omega\times [t,T]$ with respect to the measure
$dA_s(\omega)\P(d\omega)$.

This holds for instance if $U$ is a compact metric space  and the functions
$r_t(\omega,x,\cdot), l_t(\omega,x,\cdot):U\to \R$ are
continuous for every $\omega\in\Omega$, $t\in [0,T]$, $x\in
K$.
\end{hypothesis}

In this situation
Theorem \ref{Th:ex-un}  can still be applied to find a unique
solution $(Y_s^{t,x}, Z_s^{t,x})_{s\in [t,T]}$. Let us now extend
the process $Z^{t,x}$ setting $Z_s^{t,x}=0$ for $s\in [0,t)$.
The corresponding admissible control $\underline{u}^{Z,t,x}$, whose
existence is required in Hypothesis \ref{hyp:hamiltonianabis}, will be denoted
$ u^{*,t,x}$ (we set  $u^{*,t,x}(\omega,s)$
 equal
 to an arbitrary constant element of $U$ for $s\in [0,t)$).

\begin{theorem}\label{teoremacontrollononmarkovbis}
Assume that Hypotheses \ref{hyp:controllo} and
\ref{hyp:hamiltonianabis} are satisfied
  and that
\begin{equation}\label{ATintegrabiledue}
\E \exp \left((3+C_r^4) A_T\right)<\infty.
\end{equation}
Suppose also that there exists $\beta$ such that
\begin{equation}\label{compatibbetadue}
  \beta>\sup|r-1|^2, \quad \E \exp \left(\beta A_T\right)<\infty,
\quad \E[|g(X_T^{t,x})|^2e^{\beta A_T} ]<\infty, \qquad t\in
[0,T],\,x\in K,
\end{equation}
(in particular,  \eqref{compatibbetadue} follows from
\eqref{ATintegrabiledue}
 with
$\beta=3+C_r^4$ if $g$ is bounded). For any $(t,x)\in [0,T]\times K$
let $(Y^{t,x},Z^{t,x})$ denote
 the solution of the BSDE (\ref{bsdecontrollobis})
and
$u^{*,t,x}=\underline{u}^{Z,t,x}$ the corresponding admissible control.

Then  $u^*(\cdot)$ is  optimal and $Y_t^{t,x}$
is the optimal cost, i.e. $Y_t^{t,x}= J_t(x,u^*(\cdot))= v(t,x)$
$\P$-a.s.
\end{theorem}

The proof of Theorem \ref{teoremacontrollononmarkovbis} is entirely analogous
to the proof of Theorem \ref{teoremacontrollononmarkov}, the only difference
being that in the BSDE one takes the conditional expectation $\E_u^{\calf_t}$
 instead of the expectation $\E_u$.

\begin{remark}\begin{em}
\begin{enumerate}\item
Let $u\in\cala$. Then, under $\P_u$, the compensator of the
process $N$ is $A^u$ defined in (\ref{aucompensato}). It might
therefore be more natural to define as the cost corresponding to
$u\in\cala$ the  functional
$$
\E_u\left[
\int_0^Tl_t(X_t,u_t)\,dA_t^u + g(X_T)\right]=
\E_u\left[
\int_0^Tl_t(X_t,u_t)
\int_K r_t(y,u_t)\,\phi_t(dy)
\,dA_t + g(X_T)\right],
$$
instead of $J(u(\cdot))$.
This cost functional has the same form as
 $J(u(\cdot))$, with the function $l$ replaced by
 $l^0_t(x,u) :=l_t(x,u)\int_K r_t(y,u)\,\phi_t(dy)$.  Since $l^0$ is
 $\calp\otimes \calk\otimes \calu$-measurable and bounded,
 the statements of
Theorems \ref{teoremacontrollononmarkov} and \ref{teoremacontrollononmarkovbis}
remain true without any change.

\item
Suppose that the cost functional has the form
$$
J^1(u(\cdot))=\E_u\left[\sum_{n\ge1\,:\,T_n\le T}c(T_n,X_{T_n}, u_{T_n})\right],
$$
for some given  function $c:\Omega\times[0,T]\times K\times U\to
\R$ which is assumed to be bounded and  $\calp\otimes \calk\otimes
\calu$-measurable. It is well known (see e.g. \cite{B}, chapter
VII, \S 1, remark $(\beta)$) that we can reduce this control
problem to the previous one noting that
$$
J^1(u(\cdot))=\E_u \int_0^T\int_K  c(t,y, u_{t})\,p(dt\,dy)= \E_u
\int_0^T\int_K  c(t,y, u_{t})\,r_t(y,u_t)\,\phi_t(dy)\,dA_t.
$$
Thus,  $J^1(u(\cdot))$ has the same form as
 $J(u(\cdot))$, with $g=0$ and the function $l$ replaced by
 $l^1_t(x,u) := \int_K c(t,y, u)\,r_t(y,u)\,\phi_t(dy)$.  Since $l^1$ is
 $\calp\otimes \calk\otimes \calu$-measurable and bounded,
Theorems \ref{teoremacontrollononmarkov} and \ref{teoremacontrollononmarkovbis}
can still be applied.

Similar considerations obviously hold for cost functionals of the
form $J(u(\cdot))+J^1(u(\cdot))$.
\end{enumerate}
\end{em}
\end{remark}

\section{The stochastic Hamilton-Jacobi-Bellman equation}
\label{sec-HJB}

Throughout this section we still assume that a marked point
process is given, satisfying the assumptions of Section
\ref{notations}. In particular we suppose that $T_n\to\infty$
$\P$-a.s. and that (\ref{acontinuo}) holds.

We address the same optimal control problem as in the previous
section.  The associated stochastic Hamilton-Jacobi-Bellman
equation (HJB for short) is a backward stochastic differential
equation for unknown random fields on $[0,T]\times K$, having
   the Hamiltonian function   defined in (\ref{defhamiltonian}) as
  a nonlinear term. Before
  introducing the HJB equation we need a preliminary result which
  may have an interest in its  own and
 will be used to clarify the connections
with the optimal control problem and the BSDEs introduced in the
previous section, as well as in the proof of the main result,
Theorem \ref{theoremhjb}.

\subsection{A lemma of Ito type}

The Ito formula for processes defined by stochastic integrals with
respect to random measures is certainly known, see e.g.
\cite{I-W}:  it gives a canonical decomposition of $v(t,X_t)$ for
a deterministic functions $v(t,x)$ smooth enough. We need an
extension to the case when $v(t,x)$ is stochastic and itself
defined by integrals with respect to random measures. The
following result is therefore the analogue to the so-called
Ito-Kunita formula (also attributed to Bismut and Wentzell, see
e.g. \cite{Bis}, \cite{W},\cite{K}).

\begin{lemma}\label{itolemma}
Assume that $v,f:\Omega\times [0,T] \times K\to \R$ are
$Prog\otimes \calk$-measurable, $V:\Omega\times [0,T] \times
K\times K\to \R$ is $\calp\otimes \calk\otimes \calk$-measurable,
and, $\P$-a.s.
\begin{equation}\label{ipperito}
  \int_0^T| f(t,x)|\,dA_t +\int_0^T\int_K  |V(t,x,y)|\,
\phi_t(dy)\,dA_t<\infty,\qquad x\in K.
\end{equation}
Suppose that,  $\P$-a.s.
\begin{equation}\label{vVperito}
v(t,x)-v(0,x)= \int_0^t f(s,x)\,dA_s+ \int_0^t\int_K  V(s,x,y)\,
q(ds\,dy), \qquad t\in [0,T], \,x\in K.
\end{equation}
Then, $\P$-a.s.
\begin{equation}\label{itoprima}
  \begin{array}{lll}
v(t,X_t)-v(0,X_0)&\!\!=\!\!&\dis\int_0^t f(s,X_s)\,dA_s+
\int_0^t\int_K \Big(v(s-,y)-v(s-,X_{s-})+ V(s,y,y)\Big)\,
p(ds\,dy)
\\
&&\dis -\int_0^t\int_K V(s,X_s,y)\, \phi_s(dy)\,dA_s, \qquad t\in
[0,T], \,x\in K.
 \end{array}
\end{equation}

If, in addition,
$$
\int_0^T\int_K  |v(t,y)+ V(t,y,y)|\, \phi_t(dy)\,dA_t<\infty,
\qquad \P-a.s.
$$
then, $\P$-a.s.
\begin{equation}\label{itoseconda}
  \begin{array}{lll}
v(t,X_t)-v(0,X_0)&\!\!=\!\!&\dis\int_0^t f(s,X_s)\,dA_s+
\int_0^t\int_K \Big(v(s-,y)-v(s-,X_{s-})+ V(s,y,y)\Big)\,
q(ds\,dy)
\\
&&\dis +\int_0^t\int_K \Big(v(s,y)-v(s,X_{s})+ V(s,y,y)-
V(s,X_s,y)\Big)\, \phi_s(dy)\,dA_s,
 \end{array}
\end{equation}
for every $t\in [0,T]$, $x\in K$.
\end{lemma}

\begin{remark}\begin{em}
\begin{enumerate}\item
It follows from \eqref{vVperito} that $\P$-a.s. the trajectories
$v(\cdot,x)$ are cadlag for every $x\in K$. Therefore the process
$(v(t-,x))$ is well defined and $\calp\otimes \calk$-measurable.

\item
 We note that
$$
\int_0^T\int_K |V(t,X_t,y)|\, \phi_t(dy)\,dA_t= \sum_{n\ge 1}
\int_{T_{n-1}\wedge T}^{T_{n}\wedge T}\int_K |V(t,\xi_{n-1},y)|\,
\phi_t(dy)\,dA_t <\infty, \quad \P-a.s.
$$
This follows from the assumption \eqref{ipperito}, and the fact
that  the sum is finite $\P$-a.s.  due to the assumption that
$T_n\to\infty$.
 Similarly,
$$
\int_0^T |f(t,X_t)|\,dA_t + \int_0^T |v(t,X_t)|\,dA_t<\infty,
\qquad \P-a.s.
$$
so that all the integrals above are well defined: compare the
discussion in subsection \ref{subsec-stocint}.

\end{enumerate}
\end{em}
\end{remark}

\noindent{\bf Proof.} Noting that there are $N_t$ jump times $T_n$
in the time interval $[0,t]$ we have
$$
v(t,X_t)-v(0,X_0)= \sum_{n=1}^{N_t}\Big(
v(T_n-,X_{T_n})-v(T_{n-1}-,X_{T_{n-1}})\Big) + v(t,X_t) -
v(T_{N_t}-,X_{T_{N_t}}),
$$
where we use the convention $v(0-,x)=v(0,x)$. Since $X_t=
X_{T_{N_t}}$ we have
$$
v(t,X_t)-v(0,X_0)= I+II,
$$
where
$$
I= \sum_{n=1}^{N_t}\Big(
v(T_n-,X_{T_n})-v(T_{n}-,X_{T_{n-1}})\Big),
$$
$$II=
\sum_{n=1}^{N_t}\Big(
v(T_n-,X_{T_{n-1}})-v(T_{n-1}-,X_{T_{n-1}})\Big) +
v(t,X_{T_{N_t}}) - v(T_{N_t}-,X_{T_{N_t}}).
$$
Letting $H$ denote the $\calp\otimes\calk$-measurable process
$$
H_s(y)=v(s-,y)-v(s-,X_{s-}),
$$
with  the convention $X_{0-}=X_0$, we have
$$
I= \sum_{n\ge 1:T_n\le t} \Big(
v(T_n-,X_{T_n})-v(T_{n}-,X_{T_{n-1}})\Big)= \sum_{n\ge 1:T_n\le t}
H_{T_n}(X_{T_n})= \int_0^t \int_KH_s(y)\,p(ds\,dy).
$$

For $n=1,\ldots, N_t$, recalling that
$q(dt\,dy)=p(dt\,dy)-\phi_t(dy)\,dA_t$ and the definition of $p$,
$$
\begin{array}{l}
v(T_n-,x)-v(T_{n-1}-,x)
\\
\qquad\dis = V(T_{n-1},x, \xi_{n-1})- \int_{T_{n-1}}^{T_{n}}\int_K
V(s,x,y)\, \phi_s(dy)\,dA_s +\int_{T_{n-1}}^{T_{n}} f(s,x)\,dA_s
\end{array}
$$
Setting $x=X_{T_{n-1}}=\xi_{n-1}$, noting that $X_s= X_{T_{n-1}}$
for $s\in ({T_{n-1}},{T_{n}})$ and recalling that $A$ is assumed
to be continuous,
$$
\begin{array}{l}
v(T_n-,X_{T_{n-1}})-v(T_{n-1}-,X_{T_{n-1}})
\\
\qquad\dis = V(T_{n-1},\xi_{n-1}, \xi_{n-1})-
\int_{T_{n-1}}^{T_{n}}\int_K  V(s,X_s,y)\, \phi_s(dy)\,dA_s
+\int_{T_{n-1}}^{T_{n}} f(s,X_s)\,dA_s.
\end{array}
$$
Similarly,
$$
\begin{array}{l}
v(t,X_{T_{N_t}}) - v(T_{N_t}-,X_{T_{N_t}})
\\
\qquad\dis = V(T_{N_t},\xi_{N_t}, \xi_{N_t})-
\int_{T_{N_t}}^{t}\int_K  V(s,X_s,y)\, \phi_s(dy)\,dA_s
+\int_{T_{N_t}}^{t} f(s,X_s)\,dA_s.
\end{array}
$$
It follows that
$$
\begin{array}{lll}
II&=&\dis \sum_{n\ge 1:T_n\le t} V(T_{n},\xi_{n}, \xi_{n})-
\int_{0}^{t}\int_K  V(s,X_s,y)\, \phi_s(dy)\,dA_s +\int_{0}^{t}
f(s,X_s)\,dA_s
\\
&=&\dis \int_{0}^{t}\int_K  V(s,y,y)\, p(ds\,dy)
-\int_{0}^{t}\int_K  V(s,X_s,y)\, \phi_s(dy)\,dA_s +\int_{0}^{t}
f(s,X_s)\,dA_s,
\end{array}
$$
and \eqref{itoprima} is proved.  Using again the equality
$q(dt\,dy)=p(dt\,dy)-\phi_t(dy)\,dA_t$ and the additional
assumption, \eqref{itoseconda} follows as well. \qed

\begin{remark}
\begin{em}
In differential form, under the assumptions of the lemma, if
$$
dv(t,x)= f(t,x)\,dA_t+\int_K V(t,x,y)\,q(dt\, dy)
$$
then
$$
\begin{array}{lll}
dv(t,X_t) &=&\dis  f(t,X_t)\,dA_t+
 \int_K \Big(v(t-,y)-v(t-,X_{t-})+ V(t,y,y)\Big)\, q(dt\,dy)
\\
&&\dis + \int_K \Big(v(t,y)-v(t,X_{t})+ V(t,y,y)-
V(t,X_t,y)\Big)\, \phi_t(dy)\,dA_t.
 \end{array}
$$
\end{em}
\end{remark}

\subsection{The equation}

In the rest of this section we will suppose   that $U,l,r,g$ are
given satisfying Hypotheses \ref{hyp:controllo} and
\ref{hyp:hamiltoniana} as before. For technical reasons we will
also assume that the space $K$ is finite or countable (and $\calk$
is the collection of all its subsets).
 We next present
  the HJB equation by first introducing  the
  space of processes where we seek its solution.

A pair $(v,V)$ is said to belong to the space
$\H_\beta$, where $\beta\in\R$, if
\begin{enumerate}
\item $v:\Omega\times [0,T] \times K\to \R$ is $Prog\otimes \calk$-measurable,
$V:\Omega\times [0,T] \times K\times K\to \R$ is $\calp\otimes \calk\otimes \calk$-measurable;
\item
The following is finite:
$$
\begin{array}{lll}
|||(v,V)|||^2_\beta&:=&\dis \sup_{x\in K}\E\int_0^Tv(t,x)^2e^{\beta A_t}dA_t
+\E\int_0^Tv(t,X_t)^2e^{\beta A_t}dA_t
\\&&\dis
+ \sup_{x\in K}\E\int_0^T\int_KV(t,x,y)^2\phi_t(dy)\,e^{\beta A_t}dA_t
\\&&\dis
+\E\int_0^T\int_K|v(t,y)+V(t,y,y)|^2\phi_t(dy)\,e^{\beta A_t}dA_t.
\end{array}
$$
\end{enumerate}
The space $\H_\beta$, endowed with the norm $|||\cdot|||_\beta$,
is Banach space, provided we identify pairs of processes whose
difference has norm zero.

Let $f$ be the Hamiltonian function   defined in
(\ref{defhamiltonian}). A pair $(v,V)\in \H_\beta$ is called a
solution to the stochastic HJB equation if, $\P$-a.s.,
\begin{equation}\label{HJBstoch}
\begin{array}{l}
  \dis
  v(t,x)+ \int_t^T\int_KV(s,x,y)\,q(ds\, dy)
  \\\qquad\dis
  =
  g(x)
  +
  \int_t^T
\int_K\Big( v(s,y)-v(s,x) + V(s,y,y)-V(s,x,y)\Big)\,\phi_s(dy) \,dA_s
\\\qquad\dis
+\int_t^T f\Big(s,x,v(s,\cdot)-v(s,x)+ V(s,\cdot,\cdot)\Big)\,dA_s,
\qquad t\in [0,T],\, x\in K.
\end{array}
\end{equation}

We will also use the
differential notation:
$$
\left\{\begin{array}{lll}
-dv(t,x)&=&\dis
-\int_KV(t,x,y)\,q(dt\, dy)
\\&&\dis
+
\int_K\Big( v(t,y)-v(t,x) + V(t,y,y)-V(t,x,y)\Big)\,\phi_t(dy) \,dA_t
\\&&\dis
+ f\Big(t,x,v(t,\cdot)-v(t,x)+ V(t,\cdot,\cdot)\Big)\,dA_t,
\\
v(T,x)&=&g(x), \qquad t\in [0,T],\, x\in K.
\end{array}
\right.
$$

The basic result, which we assume for the moment and we will prove later, is the
following. Let $\beta_0>1$ satisfy
$$
\frac{2(2L^2+3)}{\beta_0-1}+
\frac{8(2L^2 +3)}{\beta_0}\left(1+ \frac{1}{\beta_0}\right)<1.
$$

\begin{theorem}\label{theoremhjb}
Let $K$ be finite or countable and let Hypotheses
\ref{hyp:controllo} and \ref{hyp:hamiltonianabis} be verified. Suppose that there exists $\beta$ such that
\begin{equation}\label{ipoteomain}
  \beta\ge \beta_0,\qquad
\sup_{x\in K}\E\,\left[g(x)^2e^{\beta A_T}\right]<\infty.
\end{equation}
Then the HJB equation has a unique solution
$(v,V)$ in $\H_\beta$.
\end{theorem}

\subsection{Application to control problems and BSDEs}

For  every $(t,x)\in [0,T]\times K$ we
 consider again
 the optimal control problem described just before Theorem \ref{teoremacontrollononmarkovbis}
 and  the BSDE
(\ref{bsdecontrollobis}) for the unknown processes
$(Y_s^{t,x}, Z_s^{t,x})_{s\in [t,T]}$.

Let $(v,V)\in\H_\beta$ the solution to the HJB equation
constructed in  Theorem \ref{theoremhjb}.

Then we obtain the following result.

\begin{theorem}\label{hjbcontrollo} We make the same
assumptions as in  Theorem \ref{theoremhjb},
assuming in addition that $\beta$ also satifies \eqref{compatibbetadue}.
Then for every $(t,x)\in [0,T]\times K$
we have
\begin{equation}\label{identificazioneyez}
    Y_s^{t,x}=v(s,X_s^{t,x}),
\qquad
Z_s^{t,x}(y)= v(s-,y)-  v(s-,X^{t,x}_{s-})+  V(s ,y,y).
\end{equation}
In particular, $v(t,x)=Y_t^{t,x}$ $\P$-a.s.

If  \eqref{ATintegrabiledue} also holds then
 $v(t,x)$ coincides with the
value function of the optimal control problem i.e. $
v(t,x)=\mathop{\rm ess\ inf}_{u(\cdot)\in\cala }J_t(x,u(\cdot))$
$\P$-a.s.
\end{theorem}

Equalities (\ref{identificazioneyez}) should be understood up to
sets of measure zero in $\Omega\times [t,T]$, the measure being
$dA_s(\omega)\P(d\omega)$ for the first equality and
$\phi_s(\omega,dy)dA_s(\omega)\P(d\omega)$ for the second
equality.

\noindent {\bf Proof.} We use a straightforward extension of the
Ito lemma \ref{itolemma}   to compute the stochastic differential
$dv(s,X_s^{t,x})$ on the interval $ [t,T]$ instead of $ [0,T]$.
Using the Lipschitz character of $f$ it is not difficult to check
that all the assumptions of the lemma are verified. For instance,
we check that for every $x\in K$
$$
\E\int_0^T \int_K |V(t,x,y)|\,\;\phi_t(dy)\, dA_t \le \left(
\E\int_0^T \int_K |V(t,x,y)|^2\,\;\phi_t(dy)\,e^{\beta A_t}
dA_t\right)^\frac12\! \left( \E\int_0^T e^{-\beta A_t}
dA_t\right)^ \frac12
$$ is finite, since $(v,V)\in \H_\beta$  and $\int_0^T
e^{-\beta A_t} dA_t=\beta^{-1}(1-e^{-\beta A_T})\le \beta^{-1}$,
so that $V$ satisfies the required condition \eqref{ipperito}. The
other verifications are similar and are therefore omitted.

The Ito lemma then yields
$$
\begin{array}{l}\dis
    v(s,X_s^{t,x})+\int_s^T\int_K\Big(v(r-,y)-  v(r-,X^{t,x}_{r-})+  V(r ,y,y)\,q(dr\,dy)
    \\\dis
    =
g(X_T^{t,x}) +\int_s^T f_r(X_r^{t,x},v(r-,\cdot)-  v(r-,X^{t,x}_{r-})+  V(r ,\cdot,\cdot))\,dA_r,
\qquad s\in [t,T].
\end{array}
$$
Comparing with equation (\ref{bsdecontrollobis}) and setting
$$
    \tilde Y_s^{t,x}= v(s,X_s^{t,x}),\qquad   \tilde Z_s^{t,x}
    =v(s-,y)-  v(s-,X^{t,x}_{s-})+  V(s ,y,y),
$$
 we conclude that
 the pairs
$(Y_s^{t,x},    Z_s^{t,x})$ and
$     (\tilde Y_s^{t,x},   \tilde Z_s^{t,x})$
are solutions to the same BSDE, and the latter also belongs to
$\K^\beta$ as it follows easily from the fact that $(v,V)$ belongs to
$\H_\beta$. By uniqueness for the solution to the BSDE,
(\ref{identificazioneyez}) holds.

 All the other statements follow from Theorem
\ref{teoremacontrollononmarkovbis}. \qed

\subsection{Proof of Theorem \ref{theoremhjb}}

It is convenient to first state the following simple preliminary
result.

\begin{lemma}\label{hjblineare}
Suppose
$$
-dv(t,x)=
-\int_KV(t,x,y)\,q(dt\, dy)
+
\int_KU(t,x,y)\,\phi_t(dy) \,dA_t +u(t,x)\,dA_t,
\quad
v(T,x)=g(x).
$$
Then, setting $c_\beta=\frac{2}{\beta-1} $ for $ \beta>1$, we have, for every $x\in K$,
$$
\begin{array}{l}\dis
\E\int_0^Tv(s,x)^2e^{\beta A_s}dA_s+
\E\int_0^T\int_KV(s,x,y)^2\phi_s(dy)\,e^{\beta A_s}dA_s
\\\dis
\le\E\,\left[g(x)^2e^{\beta A_T}\right]+c_\beta
 \E\int_0^Tu(s,x)^2\,e^{\beta A_s}dA_s+
 c_\beta
 \E\int_0^T\int_KU(t,x,y)^2\,\phi_s(dy)\,e^{\beta A_s}dA_s.
\end{array}
$$
\end{lemma}

\noindent{\bf Proof.} Using the identity \eqref{identitaenergia}
of Lemma \ref{BSDEcasobanale} we have
$$
\begin{array}{l}\dis
\E\,\left[v(t,x)^2e^{\beta A_t}\right]+
\beta \E\int_t^Tv(s,x)^2e^{\beta A_s}dA_s+
\E\int_t^T\int_KV(s,x,y)^2\phi_s(dy)\,e^{\beta A_s}dA_s
\\\qquad\dis
=\E\,\left[g(x)^2e^{\beta A_T}\right]+
2\E\int_t^Tv(s,x)\left[\int_KU(t,x,y)\,\phi_s(dy)+u(s,x)\right]\,e^{\beta A_s}dA_s.
\end{array}
$$
Setting $t=0$ and
using the elementary inequality
$$
2v(s,x)\left[\int_KU(t,x,y)\,\phi_s(dy)+u(s,x)\right]
\le (\beta-1)v(s,x)^2 + c_\beta \left[\int_KU(t,x,y)^2\,\phi_s(dy)+u(s,x)^2\right]
$$
the conclusion follows immediately.
\qed

\noindent{\bf Proof of Theorem \ref{theoremhjb}.}
We define a map $\Gamma:\H_\beta\to \H_\beta$ setting $(v,V)=\Gamma(u,U)$,
for $(u,U)\in \H_\beta$,
if $(v,V)$ is the solution of
$$
\left\{\begin{array}{lll}
-dv(t,x)&=&\dis
-\int_KV(t,x,y)\,q(dt\, dy)
\\&&\dis
+
\int_K\Big( u(t,y)-u(t,x) + U(t,y,y)-V(t,x,y)\Big)\,\phi_t(dy) \,dA_t
\\&&\dis
+ f\Big(t,x,u(t,\cdot)-u(t,x)+ U(t,\cdot,\cdot)\Big)\,dA_t,
\\
v(T,x)&=&g(x), \qquad t\in [0,T],\, x\in K.
\end{array}
\right.
$$
Note the two occurrences of $V$ in the right-hand side.
For fixed $x\in K$, the existence of processes $v(\cdot,x),V(\cdot,x,\cdot)$
solution to this equation follows from an application of Theorem
 \ref{Th:ex-un}. Since $K$ is assumed to be at most countable,
 the corresponding integral equation holds simultaneously
 for every $t\in [0,T]$ and $x\in K$, with the exception of
 a $\P$-null set.  The rest of the proof
consists in showing that $(v,V)\in \H_\beta$ and that $\Gamma$ is
a contraction for sufficiently large $\beta$. We limit ourselves
to proving the contraction property, since the  fact that
$(v,V)\in \H_\beta$ can be verified by similar and simpler
arguments.

Let $(u^i,U^i)\in \H_\beta$ for $i=1,2$ and let
$(v^i,V^i)=\Gamma(u^i,U^i)$. Define $\bar v=v^2-v^1$,
$\bar V=V^2-V^1$,
 $\bar u=u^2-u^1$,
$\bar U=U^2-U^1$,
$$
\bar f(t,x)= f\Big(t,x,u^2(t,\cdot)-u^2(t,x)+ U^2(t,\cdot,\cdot)\Big)
- f\Big(t,x,u^1(t,\cdot)-u^1(t,x)+ U^1(t,\cdot,\cdot)\Big).
$$
Then
\begin{equation}\label{hjbdifferenza}
    \left\{\begin{array}{lll}
-d\bar v(t,x)&=&\dis
-\int_K\bar V(t,x,y)\,q(dt\, dy)
\\&&\dis
+
\int_K\Big( \bar u(t,y)-\bar u(t,x) + \bar U(t,y,y)-\bar V(t,x,y)\Big)\,\phi_t(dy) \,dA_t
+\bar f (t,x)\,dA_t,
\\
v(T,x)&=&0, \qquad t\in [0,T],\, x\in K.
\end{array}
\right.
\end{equation}
From Lemma \ref{hjblineare} it follows that,
for every $x\in K$, $\beta >1$
$$
\begin{array}{l}\dis
\E\int_0^T\bar v(s,x)^2e^{\beta A_s}dA_s+
\E\int_0^T\int_K\bar V(s,x,y)^2\phi_s(dy)\,e^{\beta A_s}dA_s
\le
\frac{2}{\beta -1}
 \E\int_0^T\bar f (s,x)^2\,e^{\beta A_s}dA_s
\\\dis
+
 \frac{2}{\beta -1}
 \E\int_0^T\int_K\left[
 \bar u(s,y)-\bar u(s,x) + \bar U(s,y,y)-\bar V(s,x,y)\right]
 ^2\,\phi_s(dy)\,e^{\beta A_s}dA_s.
\end{array}
$$
By the Lipschitz condition on $f$ we have
\begin{equation}\label{fhamlip}
    \bar f(s,x)^2\le L^2  \int_K
\left[
 \bar u(s,y)-\bar u(s,x) + \bar U(s,y,y)\right]^2
 \,\phi_s(dy),
\end{equation}
 and it follows that,
for every $x\in K$, $\beta>1$
$$
    \begin{array}{l}\dis
\E\int_0^T\bar v(s,x)^2e^{\beta A_s}dA_s+
\E\int_0^T\int_K\bar V(s,x,y)^2\phi_s(dy)\,e^{\beta A_s}dA_s
\\\dis
\leq \frac{2(2L^2+3)}{\beta-1}\left(
 \E\int_0^T\bar u (s,x)^2\,e^{\beta A_s}dA_s
+
 \E\int_0^T\int_K\left[
 \bar u(s,y)+ \bar U(s,y,y)\right]
 ^2\,\phi_s(dy)\,e^{\beta A_s}dA_s\right)
 \\\dis
\quad  +
\frac{6}{\beta-1}
 \E\int_0^T\int_K\bar V(s,x,y)
 ^2\,\phi_s(dy)\,e^{\beta A_s}dA_s.
\end{array}
$$
Setting $c^{(1)}_\beta:=\frac{2(2L^2+3)}{\beta-1}$ it follows that
\begin{equation}\label{primastimahjb}
\sup_{x\in K}\E\int_0^T\bar v(s,x)^2e^{\beta A_s}dA_s+\sup_{x\in K}
\E\int_0^T\int_K\bar V(s,x,y)^2\phi_s(dy)\,e^{\beta A_s}dA_s
\le
 c^1_\beta\,
 |||(
 \bar u, \bar U)|||^2_\beta.
\end{equation}
We set now
$$
\bar Y_s=\bar v(s,X_s),\quad
\bar Z_s(y)=\bar v(s-,y)-\bar v(s-,X_{s-})+\bar V(s ,y,y).
$$
Recalling (\ref{hjbdifferenza}) and
applying the Ito formula of Lemma \ref{itolemma} we obtain
$$
\left\{\begin{array}{lll}
d\bar Y_t&=&\dis
\int_K\bar Z_t(y)\,q(dt\, dy)
-\bar f (t,X_t)\,dA_t
\\&&\dis
+
\int_K\Big(\bar Z_t(y)- \bar u(t,y)+\bar u(t,X_t) - \bar U(t,y,y)\Big)\,\phi_t(dy) \,dA_t,
\end{array}
\right.
$$
and $\bar Y_T=0$. Note that the term $\bar V(t,X_t,y)$
has disappeared. Using the estimate
 \eqref{stimaprlineare} in Lemma \ref{BSDEcasobanale} on the BSDE we have
$$
\begin{array}{l}\dis
\E\int_0^T\bar Y_s^2e^{\beta A_s}dA_s+
\E\int_0^T\int_K\bar Z_s(y)^2\phi_s(dy)\,e^{\beta A_s}dA_s
\le
\frac{8}{\beta}\left(1+ \frac{1}{\beta}\right)
 \E\int_0^T\bar f (s,X_s)^2\,e^{\beta A_s}dA_s
\\\dis
+
 \frac{8}{\beta}\left(1+ \frac{1}{\beta}\right)
 \E\int_0^T\int_K\left[
 \bar Z_s(y)- \bar u(s,y)+\bar u(s,X_s) - \bar U(s,y,y)
 \right]
 ^2\,\phi_s(dy)\,e^{\beta A_s}dA_s.
\end{array}
$$
Using again the inequality (\ref{fhamlip}) we obtain
\begin{eqnarray*}
\E\int_0^T\bar Y_s^2e^{\beta A_s}dA_s & + &
\E\int_0^T\int_K\bar Z_s(y)^2\phi_s(dy)\,e^{\beta A_s}dA_s \\
& \le &
\frac{8(2L^2 +3)}{\beta}\left(1+ \frac{1}{\beta}\right)
 \E\int_0^T\bar u (s,X_s)^2\,e^{\beta A_s}dA_s
\\
&+&
 \frac{8(2L^2 +3)}{\beta}\left(1+ \frac{1}{\beta}\right)
 \E\int_0^T\int_K\left[
\bar u(s,y) + \bar U(s,y,y)
 \right]
 ^2\,\phi_s(dy)\,e^{\beta A_s}dA_s
 \\
 &+ &  \frac{24}{\beta}\left(1+ \frac{1}{\beta}\right)
 \E\int_0^T\int_K
 \bar Z_s(y)
 ^2\,\phi_s(dy)\,e^{\beta A_s}dA_s.
\end{eqnarray*}

Setting $c^{(2)}_\beta:=\frac{8(2L^2 +3)}{\beta}\left(1+ \frac{1}{\beta}\right)$ it follows that
\begin{equation}\label{ybarzbar}
\E\int_0^T\bar Y_s^2e^{\beta A_s}dA_s+
\E\int_0^T\int_K\bar Z_s(y)^2\phi_s(dy)\,e^{\beta A_s}dA_s
\le
c^2_\beta \,
 |||(
 \bar u, \bar U)|||^2_\beta.
\end{equation}
Recalling the definition of
$\bar Y,\bar Z$ and using the fact that $A$ is assumed
to be continuous we have
\begin{equation}\label{ybarzbardue}
\begin{array}{l}\dis
 \E\int_0^T\int_K\left[
\bar v(s,y) + \bar V(s,y,y)
 \right]
 ^2\,\phi_s(dy)\,e^{\beta A_s}dA_s
 =
\E\int_0^T\int_K
\left[
\bar Z_s(y)+\bar Y_s\right]^2\phi_s(dy)\,e^{\beta A_s}dA_s
\\\dis
\le
2\E\int_0^T\bar Y_s^2e^{\beta A_s}dA_s+
2\E\int_0^T\int_K\bar Z_s(y)^2\phi_s(dy)\,e^{\beta A_s}dA_s
\le
c^2_\beta \,
 |||(
 \bar u, \bar U)|||^2_\beta,
\end{array}
\end{equation}
where the last inequality is due to (\ref{ybarzbar}).
Recalling that
$\bar Y_s=\bar v(s,X_s)$,
it follows from
(\ref{primastimahjb}),   (\ref{ybarzbar}), (\ref{ybarzbardue}) that
$ |||(\bar v, \bar V)|||^2_\beta\le
c_\beta \,
 |||( \bar u, \bar U)|||^2_\beta$ where $c_\beta= c^{(1)}_\beta+ c^{(2)}_\beta$ is
 $<1$ by the assumptions.
 This proves the required contraction property
 and finishes the proof.
\qed

\bigskip

{\bf Acknowledgements.} We wish to thank Prof. Ying Hu for
pointing out to us the reference \cite{Xia} and Prof. Huy\^{e}n Pham
for discussions on connections between BSDEs and point processes.

\end{document}